\documentclass[10pt,reqno,oneside]{amsart}
\usepackage{xcolor} 
\usepackage{amssymb,amsmath,mathrsfs,amsthm}
\usepackage{enumitem}
\usepackage{mathtools}
\usepackage{geometry}
\usepackage{marginnote}
\usepackage{cancel}
\usepackage{amsfonts}
\allowdisplaybreaks
  \chardef\forshowkeys=0
  \chardef\refcheck=0
  \chardef\showllabel=0
\ifnum\forshowkeys=1
   \usepackage[notref,notcite,color]{showkeys}
\fi
\usepackage[unicode,breaklinks=true,colorlinks=true,linkcolor=blue,urlcolor=blue,citecolor=blue]{hyperref}

\ifnum\showllabel=1
\def\llabel#1{\marginnote{\color{gray}\rm(#1)}[-0.0cm]\notag}
\else
\def\llabel#1{\notag}
\fi

\ifnum\refcheck=1
\usepackage{yfonts}\usepackage{refcheck}
\fi

\usepackage{esint,comment}

\usepackage[dvips]{graphicx}
\allowdisplaybreaks
\newtheorem{thm}{Theorem}

 \newtheorem{lemma}[thm]{Lemma}

\theoremstyle{definition}

\numberwithin{equation}{section}

\def\eqL{\stackrel{L}{=}}
\def\leqL{\stackrel{L}{\leq}}
\def\Tmax{T_\text{max}}
\def\inon#1{\hbox{\ \ \ \ \ \ \ }\hbox{#1}}

\def\indeq{\quad{}} 

\def\comma{ {\rm ,\qquad{}} }            
\def\commaone{ {\rm ,\quad{}} }          

\definecolor{colorik}{rgb}{0.35,0.1,0.6}\def\cole{\color{coloroooo}}

\definecolor{coloroooo}{rgb}{0.45,0.0,0.0}
\def\cole{\color{coloroooo}}
\def\cole{\colb}

\def\colb{\color{black}}

\def \no#1#2#3 {{\bf #1} (#3), #2.}
\def \eds#1#2#3 {#1, #2, #3.}

\def\:{{\colon}}

\def\be#1{\begin{equation}\label{#1}}
\def\ee{\end{equation}}

\def\<{\langle}
\def\>{\rangle}
\def\coloneqq{:=}

\newcommand{\na}{\nabla}

\newcommand{\lec}{\lesssim}

\newcommand{\bs}{\begin{split}}
\newcommand{\essss}{\end{split}}

\renewcommand{\div}{\operatorname{div}}
 
\newcommand{\eqnb}{\begin{equation}}
\newcommand{\eqne}{\end{equation}}

\renewcommand{\ee}{\mathrm{e}}

\newcommand{\p}{\partial}

\renewcommand{\d}{\,\mathrm{d}}

\newcommand{\supp}{\operatorname{supp}}

\newcommand{\curl} {\mathop{\mathrm{curl}}}

\begin{document}
\title[Constructing free-boundary flow with limited regularity]{Construction of the free-boundary 3D incompressible Euler flow under limited regularity}
\author[M.~Aydin]{Mustafa Sencer Aydin}
\address{Department of Mathematics, University of Southern California, Los Angeles, CA 90089}
\email{maydin@usc.edu}

\author[I.~Kukavica]{Igor Kukavica}
\address{Department of Mathematics\\
University of Southern California\\
Los Angeles, CA 90089}
\email{kukavica@usc.edu}

\author[W.O.~O\.za\'nski]{Wojciech S.~O\.za\'nski}
\address{Department of Mathematics\\Florida State University\\Tallahassee, FL 32306}
\email{wozanski@fsu.edu}

\author[A.~Tuffaha]{Amjad Tuffaha}
\address{Department of Mathematics and Statistics, American University of Sharjah, Sharjah, UAE}
\email{atufaha\char'100aus.edu}

\date{\today}
\maketitle

\begin{abstract}
We consider the three-dimensional Euler equations in a domain with a free boundary with no surface tension. We construct unique local-in-time solutions in the Lagrangian setting for $u_0 \in H^{2.5+\delta }$ such that the Rayleigh-Taylor condition holds and $\mathrm{curl}\,u_0 \in H^{2+\delta }$ in an arbitrarily small neighborhood of the free boundary. We show that the result is optimal in the sense that $H^{3+\delta}$ regularity of the Lagrangian deformation near the free boundary can be ensured if and only if initial vorticity has $H^{2+\delta}$ regularity of vorticity near the free boundary.
\end{abstract}

\section{Introduction}\label{sec_int}
We address the local-in-time well-posedness of the three-dimensional incompressible Lagrangian Euler equations
  \begin{align}
  \begin{split}
   &\partial_t v_i + a_{ki} \partial_k q = 0, \quad i = 1, 2, 3
   \\&
   a_{ik} \partial_i v_k = 0
   ,
  \end{split}
  \label{EQ75}
  \end{align}
where $\eta_t = v$ and $a = (\nabla \eta)^{-1}$, with the initial conditions given by $(v, a, \eta)(0) = (v_0, I, x)$. It is well-known that the  above system represents the Euler equations
  \begin{align}
  \begin{split}
  &u_t + u \cdot \nabla u + \nabla p = 0
  \\&
  \nabla \cdot u = 0
  \end{split}
  \llabel{EQ01}
  \end{align}
in an evolving domain $\Omega(t) \subseteq \mathbb{R}^3$ set in the Lagrangian variable $\eta$, where $\eta_t(x, t) = u(\eta(x, t), t)$ with $\eta(x, 0) = x$; see \cite[Section~2.1]{KO} for example.

The local existence problem for the Euler equations with an evolving boundary has a rich history. Initially, the problem of local existence was considered in \cite{Sh, Shn, Y1, Y2} under the irrotationality assumption. The need for the Rayleigh-Taylor condition on the pressure was demonstrated by Ebin in \cite{E}, who showed ill-posedness without this condition.  From a physical point of view, this condition is a requirement that the pressure in the fluid at the free boundary is higher than the pressure of the air.  In \cite{W1, W2}, S.~Wu proved the local existence for rotational Sobolev smooth initial data. By incorporating the Rayleigh-Taylor condition, \cite{ChL} provided a~priori bounds for the existence of solutions in $H^r$ for $r>3$.  In \cite{CS1} and \cite{CS2}, Coutand and Shkoller derived new energy estimates and constructed solutions to the initial value problem for $u_0\in H^4$ with $\curl\,u_0\in H^{3.5}$, see \cite[Theorem~8.4]{CS2}.  The local existence for the case of $H^3$ initial data was established independently in~\cite{CS1, SZ1, SZ2, ZZ}.  The works \cite{SZ1, SZ2, AM1, AM2} considered the limit of vanishing surface tension.  For other works on different aspects of the free-boundary fluid problems, see~\cite{ABZ2,B,CL,ChL,EL,I,KT1,KT2,KT3,KT4,L,Lin1,Lin2,MC,MR,N,T}, for the problem with surface tension, we refer the reader to~\cite{ABZ1,DE1,DE2,IK,OT,P,S}, while for the global existence results, see~\cite{AD,GMS,HIT,IP,IT,W3}.

In this work we are concerned with lowering the necessary regularity from $H^{3}$ to the $H^{2.5+}$ level, which is minimal from the viewpoint of the local well-posedness theory for $3$D Euler, even without free boundaries~\cite{BL1,BL2}. Such improvement is not only a matter of introducing fractional calculus. In fact, the main  obstacle is due to the nature of the free-boundary problem. Namely, in most of the works, the problem of constructing free-boundary flows relies heavily on an  additional regularity of the Lagrangian map, which in turn depends on the additional regularity of the vorticity. This can be seen from  the Cauchy invariance \eqref{EQ23},
\[
    \nabla ((\curl \eta)_i)
    =
    \epsilon_{ijk} (\delta_{km} - \partial_k \eta_m)\partial_j \nabla (\eta_m) 
    + t  \nabla (\omega_0)_i
    +
   2 \int_0^t
     \epsilon_{ijk} \partial_k v_m \partial_j \nabla \eta_m
        \,\d s,
  \]
which was introduced in this context in~\cite{KTV}; see also \cite{DK, KTV, KTVW, KO}. Thus, provided $\nabla v$ remains bounded, this identity shows that $\eta \in H^{3+\delta}$ \emph{if and only if} $\omega_0 \in H^{2+\delta}$. In particular, it is thanks to the initial regularity of the vorticity that we expect the Lagrangian trajectories to remain more regular, here obtaining $0.5$ Sobolev regularity more than the $2.5+\delta$ Sobolev regularity which can be expected from the definition of the trajectories for $v\in H^{2.5+\delta}$.

We also note that, in the Eulerian setting, the local existence in $H^{2.5+\delta}$ was recently obtained in \cite{WZZZ}. However, due to the lack of the chain rule in $H^{2.5+\delta}$, the results in \cite{WZZZ} do not apply to the Lagrangian setting and also do not guarantee additional regularity of the Lagrangian particle map.

Recently, in \cite{KO}, two of the authors obtained a priori bounds for solutions of \eqref{EQ75} under the minimal $H^{2.5+\delta}$ regularity on the velocity, see Lemma~\ref{L01} below.  The estimates require $H^{2+\delta}$ regularity of $\omega_0$ in any neighborhood of the free boundary. The difficulty in obtaining such a~priori bounds is that enhanced regularity is not well propagated by an evolution equation. However, in the case of the Euler equations, this is possible due to the local nature of the Cauchy invariance \eqref{EQ23}. In fact, \cite{KO} demonstrated that its localized version can guarantee extra $0.5$ Sobolev regularity of the Lagrangian trajectories $\eta$ near the free boundary \cite[(3.9)]{KO}, which in turn enables a localized a~priori estimate. 

However, we emphasize that local-in-time existence of solutions does not simply follow from the a~priori estimates. Indeed, any regularization procedure that could lead to an existence result must respect many of the symmetries and cancellations used in the obtained a priori estimates. In particular, the special structure of the integral which gives the enhanced regularity of the free boundary (see~\eqref{EQ71}) has to be preserved. This was well-addressed in \cite{CS1}, where the authors introduced the tangential mollification procedure. However, this appears to be insufficient in our approach, which additionally requires the preservation of the Cauchy invariance property. For this reason, we introduce a completely new approach: We first approximate the initial velocity field $v_0\in H^{2.5+\delta}$ by a more regular velocity field $v_0^{(r)}$, namely $v_0^{(r)}\in H^4$. We then use the existing theory to obtain local solutions for each $r>0$, and take $r\to 0^+$. The main difficulty in such approach is that the existence interval $[0,\Tmax)$ for each $r>0$ may converge to $0$ as $r\to 0^+$. In order to describe our solution to this problem, we first introduce some notation and state the main result.

\subsection{Euler Equations in the Lagrangian Setting}

The Lagrangian velocity and pressure are denoted by $v= (v_1,v_2,v_3)$ and $q(x,t)$ respectively with the Lagrangian variable~$x$. The Euler equations then become
  \begin{align}
  \begin{split}
     \p_t v_i &= - a_{ki} \p_k q,\qquad i=1,2,3,\\
     a_{ik} \p_i v_k &=0
  \end{split}
  \label{EQ03}
  \end{align}
in $\Omega \times (0,T)$, where $\Omega \coloneqq \Omega (0) = \mathbb{T}^2 \times (0,1)$. 
Here, $a_{ij}$ denotes the $(i,j)$-th entry 
of the matrix $a=(\na \eta )^{-1}$, where $\eta $ stands for the particle map, i.e., the solution of the system
  \begin{align}
  \begin{split}
     \eta_t (x,t) &= v(x,t) \\
     \eta (x,0) &= x
  \end{split} 
  \label{EQ04}
  \end{align}
in $\Omega \times [0,T)$. 
The incompressibility condition in \eqref{EQ03} implies $\mathrm{det}\, \na \eta =1$ for all times, which means that $a$ is the corresponding cofactor matrix. Therefore,
  \begin{align}
     a_{ij} =\frac12 \epsilon_{imn} \epsilon_{jkl} \partial_{m}\eta_k\partial_{n}\eta_l
  ,
  \label{EQ05}
  \end{align}
where $\epsilon_{ijk}$ denotes the permutation symbol.
The initial condition for $v$ is the same as for $u$, i.e.,
$v_0 \coloneqq v(0)=u_0$.
As for the boundary conditions, the Eulerian velocity vanishing on the bottom boundary
translates to
  \begin{align}
     v_3 =0 \qquad \text{ on }\Gamma_0,
  \label{EQ06}
  \end{align}
and the zero surface tension condition at the top boundary $\Gamma_1 \coloneqq \Gamma_1(0) = \{ x_3 =1 \}$ becomes
  \begin{align}
   q=0 \qquad \text{ on } \Gamma_1 \times (0,T).
  \label{EQ07}
  \end{align}
Observe that in the Lagrangian setting, neither $\Gamma_0$ nor $\Gamma_1$ depend on time. Finally, we denote by $\chi = \chi(x_3)$ and $\psi = \psi(x_3)$, both in $C^{\infty}(\mathbb{R},[0,1])$, smooth cut-off functions such that  $\chi$ is supported in
$\Omega'$, which is 
a fixed neighborhood of $\Gamma_1$, and $\chi = 1$ in a smaller
neighborhood. We require that $\supp \psi \Subset \{\chi = 1\}$ with
$\psi = 1$ in a neighborhood of~$\Gamma_1$.

For the rest of the paper, we abbreviate $L^p((0,T);X(\Omega))$ and
$C([0,T];X(\Omega))$ by  $L^p_TX$ and $C_TX$, respectively. We also use the symbol ``$\lec $'' to denote ``$\leq C $'', where $C>0$ is an absolute constant. Unless specified otherwise, the spatial and time domains are assumed to be $\Omega$ and $[0,T]$, respectively.

\subsection{The main result}

For the remainder of the paper, fix any $\delta\in (0,1/2]$.
We consider $v_0 \in H^{2.5+\delta }$ such that $\mathrm{curl}\,v_0\in H^{2+\delta} (\Omega')$ and such that the corresponding pressure function $q_0$, defined as the unique solution of $\Delta q_0 =- \p_i(v_0)_j \p_j (v_0)_i$ in $\Omega$ with boundary conditions $q_0=0$ on $\Gamma_1$ and $\p_3 q_0 =0 $ on $\Gamma_0$ and periodic boundary conditions in $x'$, satisfies the Rayleigh-Taylor condition~\eqref{EQ20}. 

\cole
\begin{thm}[Main result]
\label{T01}
Given $v_0$ as above such that
$
   \partial_{3} q(0) \leq -b<0
$
there exists\\
$T=T(\Vert v_0\Vert_{{2.5+\delta}}, \Vert \chi \curl v_0\Vert_{{2+\delta}}, b)>0$ and a unique solution $(v,q,\eta )$
on (0,T] to the free-boundary Euler equations \eqref{EQ03}--\eqref{EQ07} such that
  \eqnb\label{EQ08}
   \sup_{t\in [0,T] }
      \left(
         \| v \|_{2.5+\delta } +\|q \|_{2.5+\delta }
    + \| \chi q \|_{3+\delta }
    + \| \chi \eta \|_{3+\delta } \right)
   \leq
   C (\| v_0 \|_{2.5+\delta } , \| \chi \curl v_0 \|_{2+\delta },b )
  \eqne
and
  \begin{equation}
   \partial_{3} q(t) \leq -\frac{b}{2}<0
   \comma x \in \Gamma_1
   \commaone t\in(0,T)   
   .
   \label{EQ20}
  \end{equation}
Moreover, $\chi \eta \in C_T H^{s+0.5}$ and  $ \eta ,v \in  C_T H^{s}$ for every $s<2.5+\delta$.
\end{thm}
\colb

As mentioned above, such existence result is possible thanks to the a~priori bounds
from~\cite[Theorem~1~and~(3.29)]{KO}, which may be stated as follows.

\cole
\begin{lemma}[a~priori bounds]
\label{L01}
Given
$\varepsilon\in(0,1]$ and
$v_0$ as in Theorem~\ref{T01}, there exists
  \begin{equation}
   T_0= T_0 (\varepsilon , \Vert v_0\Vert_{{2.5+\delta}}, \Vert \chi \omega_0 \Vert_{{2+\delta}}, b)>0
   \llabel{EQ36}
  \end{equation}
with the following property: If $(v,q,a,\eta )\in H^4 \times H^{4.5} \times H^{3.5} \times H^{4.5}$ is a solution of the Euler system in the Lagrangian setting \eqref{EQ03}--\eqref{EQ07} on $[0,T_0]$, then
 \begin{align}
  \label{EQ15_RT}
  \p_3 q (x,t) \leq - \frac{b}2 
  \quad \text{ for } x\in \Gamma_1,  t\in [0,T_0 ] 
  \end{align}
  and
\eqnb\label{EQ13}
  \| I -a  \|_{1.5+\delta },\| I -aa^T  \|_{1.5+\delta }, \| I-\na \eta \|_{1.5+\delta }
  \leq \varepsilon,
  \eqne
  \eqnb \label{EQ14}
  \Vert \eta\Vert_{2.5+\delta}, \Vert a\Vert_{1.5+\delta} \lec 1 \qquad  \text{ for } t\in [0,T_0 ]
  ,
  \eqne
 \eqnb\label{EQ15}
 \| v(t) \|_{2.5+\delta },\| q(t) \|_{2.5+\delta }, \| \chi \eta (t) \|_{3+\delta } \lec \int_0^t P \d s + t \| \chi \omega_0 \|_{2+\delta } +\| \psi v_0 \|_{2.5+\delta }^2 + \| v_0 \|_{2+\delta }^2 + \| v_0 \|_{L^2} +1
  ,
 \eqne
for $t\in [0,T_0]$, where $P$ denotes a polynomial in $\| v \|_{2.5+\delta }$, $\| \eta \|_{2.5+\delta }$, and $\| \chi \eta \|_{3+\delta }$.
\end{lemma}
\colb

We note that the above lemma is proven in \cite{KO} under the assumption of $C^\infty$ smoothness of the solution, but it generalizes directly to the regularity assumed by the lemma. We also recall from \cite[(3.21)]{KO} that
\eqnb\label{psiqt}
  \| \psi q \|_{3+\delta },
  \| \psi q_t \|_{2.5+\delta } \leq P,
\eqne
for all $t\in [0,T_0]$.
For simplicity, we assume that $T_0\leq 1$.

\subsection{Sketch of the proof}

In order to achieve Theorem~\ref{T01},  we consider a regularization $v_0^{(r)}$ of $v_0$ (see Step~1 below), which gives rise to a unique solution on some time interval $[0,\Tmax )$, depending on $r$.
We then observe that the a~priori can be ``unlocalized'' for all $r>0$ and the entire interval $[0,T_0]$. To be more precise, if $\omega_0 \in H^2$ (rather than merely in a neighborhood of $\Gamma_1$) then one can remove the cutoff $\chi$ in \eqref{EQ15}. The point here is that, since $T_0$ depends only on $\| v_0 \|_{2.5+\delta }$ and $\| \chi \omega_0 \|_{2+\delta}$ (and not $\| \omega_0\|_{2+\delta}$), we obtain $H^{3+\delta }$ control of $\eta$.  The control will become worse as $r\to 0^+$ (via the term $t\| \omega_0 \|_{2+\delta }$, see \eqref{EQ15a} below), but it remains valid over the entire time interval $[0,T_0]$; see Step~3a below for details.

Thanks to this, we can then obtain (in Step~3b) a \emph{linear} ODE-type control of $\| \eta \|_{4.5}+\| v \|_{4}$ for all $t\in [0,\min (\Tmax , T_0 ))$ (see \eqref{EQ74b} below), that is as long as the solution evolving from the regularized data exists. Thanks to the linearity, a simple Gronwall-type estimate guarantees that $\Tmax \geq T_0$, which ensures that we can take the limit $r\to 0^+$ in function spaces over the same time interval $[0,T_0]$ (see Step~4). 

We emphasize that the implicit constant in the linear inequality might depend on the size of $\| \eta \|_{3+\delta }$, and so it might blow-up as $r\to 0^+$, but the time interval on which it holds will not shrink. We only take the limit $r\to \infty$ in the limited-regularity setting of the a~priori bound given by Lemma~\ref{L01} above.

We note that the above issue of the growing constant in the linear inequality brings another difficulty. Namely,  $T_0 $ in Lemma~\ref{L01} depends on $\varepsilon$, and one needs to be particularly careful about any absorption type tricks involving smallness \eqref{EQ13} of the deformation, which is addressed  in Section~\ref{sec_step3b}.

Another difficulty lies in the fact that the proof of the linear estimate involves some manipulations of tangential derivatives of order $5$ of both $\eta$ and $v$. These are not well-defined for $\eta\in H^{4.5}$ and $v\in H^4$ and, in order to get around this issue, one needs to perform the estimates on the level of tangential difference quotients, and a limit in the resulting estimate. This is achieved in Sections~\ref{sec_tan}~and~\ref{sec_com}.
One of the challenges in such approach is that  the product rule for the difference quotients involves the translation operator, for which we need to show that the Rayleigh-Taylor type cancellation property is not lost; see \eqref{EQ113} for details.

Finally, we prove the uniqueness part of Theorem~\ref{T01} in Section~\ref{sec_uniqueness}.

\section{Proof of Theorem~\ref{T01}: existence}\label{sec_main}

In this section, we prove the existence part of Theorem~\ref{T01}. We allow all constants to depend on
$ \Vert \bar v_0\Vert_{{2.5+\delta}}$, $ \Vert \chi \curl\bar v_0\Vert_{{2+\delta}}$, and~$b$. \\

\noindent\texttt{Step~1.}

We regularize the initial datum.\\

To this end, we use a modification of the approximation scheme from~\cite{CS2}.
First, denote the initial velocity by $\bar{v}_0$ instead of~$v_0$.
Then, for $r\in(0,1]$, let $\phi_r \in C_0^\infty (\mathbb{R}^3)$ be a family of standard mollifiers with $\supp \phi_r \subseteq B(0,r)$.
We construct the approximation by defining
  \begin{align}
  \begin{split}
   \tilde v_0^{(r)}\coloneqq\begin{pmatrix} v_0 (x_1,x_2,(x_3+r)/(1+2r)) \\
   v_0 (x_1,x_2,(x_3+r)/(1+2r)) \\
   (1+2r)^{-1}v_0 (x_1,x_2,(x_3+r)/(1+2r)) 
   \end{pmatrix}, \quad
   w^{(r)}(x) \coloneqq
   \phi_r * \tilde v_0^{(r)}
   , \quad \text{ and } \quad v_0^{(r)} \coloneqq w^{(r)} - \nabla h^{(r)},
  \end{split}
  \llabel{EQ88}
  \end{align}
where $h^{(r)}$ is the solution to the elliptic problem 
  \begin{align}
  \begin{split}
   \Delta h^{(r)}  &= \div w^{(r)}
      \inon{in $\Omega$},
   \\
   h^{(r)} &=  0
      \hspace{0.97cm}\inon{on $\Gamma_1$},
   \\
     \partial_3 h^{(r)} &= w_3^{(r)}
\hspace{0.52cm}      \inon{on $\Gamma_0$}     
  \end{split}
  \label{EQ09}
  \end{align}
with the periodic boundary conditions in  $x_1$ and~$x_2$.
Note that $\tilde v_0$ is defined and smooth in
$\mathbb{T}^2\times(-r,1+r)$,
and thus $w^{(r)}$ is smooth in a neighborhood of $\Omega$.
As $r\to 0$, the solution of \eqref{EQ09} converges to the solution of
  \begin{align}
  \begin{split}
   \Delta h
      &= 0
      \inon{in $\Omega$}
   \\
   h &=  0
      \inon{on $\Gamma_1$}
   \\
     \partial_3 h &= 0
      \inon{on $\Gamma_0$}
   ,
  \end{split}
  \label{EQ10}
  \end{align}
whence $h^{(r)}\to 0$ as $r\to0$.
By the elliptic regularity  \cite[Chapter~2]{LM}, we have $v_0^{(r)} \in H^s(\Omega)$ for $s \ge 2.5+\delta$ and $v_0^{(r)} \rightarrow \bar{v}_0$ in $H^{2.5+\delta}$ as $r \rightarrow 0$. We note that the point of the definition of $\tilde v_0^{(r)}$ is that 
  \begin{align}
  \begin{split}
   &
   \curl \tilde v_0^{(r)}
   =
   (1+2r)^{-1}
   \omega_0(x_1,x_2,(x_3+r)/(1+2r))
   ,
  \end{split}
   \label{EQ58}
  \end{align}
which gives
  \begin{align}
   \Vert \curl v_0^{(r)}\Vert_{H^{2+\delta}(\Omega')}
   = \Vert \curl w^{(r)}\Vert_{H^{2+\delta}(\Omega')}
   \lec \Vert \curl \tilde v_0^{(r)}\Vert_{H^{2+\delta}(\Omega')}
   \lec \Vert \curl v_0\Vert_{H^{2+\delta}(\Omega)}
   \lec 1.
  \label{EQ11}
  \end{align}\\

\noindent\texttt{Step~2.} For each $r>0$, we obtain a regular solution for some $\Tmax>0$.\\

Namely, for each $r>0$, we use \cite[Theorem~2.8]{CS2} to obtain $T>0$ and a solution $(v,q,\eta)$ to the Euler system \eqref{EQ03}--\eqref{EQ07} such that 
  \begin{align}
  \begin{split}
   &\eta  \in  L^\infty_T H^{4.5} \cap C_T H^{4}
   \\&
   v  \in L^\infty_T H^{4} \cap C_T H^{3.5}
   \\&
   v_t  \in L^\infty_T H^{3.5}
   \\&
   q \in L^\infty_T H^{4.5}
   \\&
   q_t  \in L^\infty_T H^{4}
   ,
  \end{split}
  \llabel{EQ12}
  \end{align}
with~\eqref{EQ20}.
We denote by $\Tmax$ the maximal time of existence of this solution. Note that, for brevity, we use the notation $(v,q,\eta)$, while this is in fact the solution obtained with the regularized initial datum $v_0^{(r)}$.
  
  The main point here is that  $\Tmax$ depends on $r$, and it seems possible that $\Tmax \to 0$ as $r\to 0^+$. In the next step, we show not only that this does not happen, but also that $\Tmax $ can be bounded below by $T_0$, which is given by the a priori estimate (Lemma~\ref{L01}). 
  
We emphasize that, for each sufficiently small $\varepsilon >0$, the a~priori estimate \eqref{EQ15} implies that
\eqnb\label{EQ15_cons}
\| v \|_{2.5+\delta } + \| q \|_{2.5+\delta } + \| \eta \|_{2.5+\delta } + \| \chi \eta \|_{3+\delta } \lec 1
,
\eqne
for all times $t\in [0,T_0]\cap [0,\Tmax)$, and for each~$r$. \\

\noindent\texttt{Step~3.} We show that there exists a sufficiently small $\varepsilon \in (0,1)$ such that $\Tmax\geq T_0$ for each $r>0$, where $T_0$ is given by Lemma~\ref{L01}. 
(Note that $T_0$ does not depend on the regularization parameter $r$, but on $\varepsilon$; recall~\eqref{EQ13}.)  \\

 We first extend the a~priori inequality \eqref{EQ15} to the case when $\omega_0 \in H^{2+\delta }$ (rather than merely $\chi \omega_0 \in H^{2+\delta }$).\\

\noindent\texttt{Step~3a.} We show in Section~\ref{sec_step3a}  that, for each $\varepsilon >0$, if the assumptions of Lemma~\ref{L01} hold and also $\| \omega_0 \|_{2+\delta }<\infty$, then 
 \eqnb\label{EQ15a}
   \|  \eta (t) \|_{3+\delta } \lec \int_0^t P \d s + t \| \omega_0 \|_{2+\delta } +\|  v_0 \|_{2.5+\delta }^2  +1
 \eqne
and 
 \eqnb\label{EQ15b}
 \| q (t) \|_{3+\delta } \lec \| v (t) \|_{2.5+\delta}^2 
 ,
 \eqne
for $t\leq T_0 $, where $P$ denotes a polynomial in~$\|  \eta \|_{3+\delta }$. (Recall from Lemma~\ref{L01} that $T_0$ depends on $\varepsilon$.)\\

Note that \eqref{EQ15a} differs from the a~priori estimate \eqref{EQ15} in that cutoffs $\chi$ and $\psi$ are replaced by~$1$. It turns out that these improved a~priori estimates \eqref{EQ15a}--\eqref{EQ15b} suffice to obtain a linear estimate on the higher order norms for local-in-time solutions from  Step~2.

We emphasize that inequalities \eqref{EQ15a}--\eqref{EQ15b} remain valid until $t=T_0$, which is independent of $r$, even though the upper bounds that they provide, including the implicit constant in \eqref{EQ15b}, could deteriorate as $r\to 0^+$.\\ 

\noindent\texttt{Step~3b.} We show in Section~\ref{sec_step3b} that if $\varepsilon >0$ is sufficiently small, then 
 \begin{align}
     \Vert \eta (t) \Vert_{{4.5}}^2 + \Vert v (t) \Vert_{{4}}^2
      \lec_r
       C_{ \Vert \eta (0) \Vert_{{4.5}}, \Vert v (0) \Vert_{{4}}}
        + \int_0^t  \left( \Vert \eta (s) \Vert_{{4.5}} + \Vert v (s) \Vert_{{4}}\right)^2 \,\d s
  \label{EQ74b}
  \end{align}
holds for $t\in [0,T_0]$, where the implicit constant may depend on $\sup_{t\in [0,T_0]} ( \| v \|_{2.5+\delta} + \| \eta \|_{3+\delta } + \| q \|_{3+\delta })$, and so in particular on~$r$.

Thanks to Step~3b, we see that, for each $r>0$, the sum $\Vert \eta(t)\Vert_{{4.5}}^2 + \Vert v (t) \Vert_{{4}}^2$  satisfies a linear integral inequality. Thus a simple Gronwall argument shows that the solution $(v,q,\eta)$ from Step~2 cannot blowup for as long as $\| v \|_{2.5+\delta}$, $\| \eta \|_{3+\delta }$ , and $\| q \|_{3+\delta }$  remain bounded. The improved a priori estimate from Step~3a shows that for each $r$ these remain bounded on $[0,T_0]$, although the bound can become larger as $r\to 0^+$. Consequently, the higher order norms
$\| \eta \|_{4.5}$, $\|v\|_4$, $\| v_t \|_{3.5}$, $\| q \|_{4.5}$, and $\| q_t \|_4$ can become larger as $r\to 0^+$, but they do not blow-up at any $t\in [0,T_0]$. As in \cite{KO}, one obtains an inequality
  \eqnb\llabel{EQ74c}
  \| q \|_{4.5},\| q_t \|_{4}, \| v_t \|_{3.5} \lec C_{\| \eta \|_{3+\delta }, \| q \|_{3+\delta }} (\| \eta \|_{4.5} + \| v \|_4)
  ,
  \eqne
for each $t\in [0,T_0]\cap [0,\Tmax)$, using the Poisson equations for $q$ and $q_t$ and the Euler equation for~$v_t$.\\

\noindent\texttt{Step~4.} We take the limit $r\to 0^+$ to obtain existence until time~$T_0 $.\\

To this end, using Step~3 and the a~priori estimate from Lemma~\ref{L01}, we can find a subsequence $r_k\to 0^+$, and relabel it back to $r$, such that 
  \begin{align}
  \begin{split}
     &
     v^{(r)} \rightharpoonup v \text{ weakly-* in }   L^\infty_{T_0} H^{2.5+\delta}
     ,
     \\&
     v^{(r)}_t \rightharpoonup v_t \text{ weakly-* in }   L^\infty_{T_0} H^{1.5+\delta}
     ,
     \\&
     \psi q^{(r)} \rightharpoonup \psi q \text{ weakly-* in }   L^\infty_{T_0} H^{3+\delta}
     ,
     \\&
     q^{(r)} \rightharpoonup q \text{ weakly-* in }   L^\infty_{T_0} H^{2.5+\delta}
     ,
     \\&
     \chi \eta^{(r)} \rightharpoonup \chi \eta \text{ weakly-* in }   L^\infty_{T_0} H^{3+\delta}
     ,
     \\&
     \eta^{(r)} \rightharpoonup \eta \text{ weakly-* in }   L^\infty_{T_0} H^{2.5+\delta}
     ,
     \\&
     \eta^{(r)}_t \rightharpoonup \eta_t \text{ weakly-* in }   L^\infty_{T_0} H^{2.5+\delta}
     ,
     \\&
     \chi^2 a^{(r)} \rightharpoonup \chi^2 a \text{ weakly-* in }   L^\infty_{T_0} H^{2+\delta}
     ,
     \\&
     a^{(r)} \rightharpoonup a \text{ weakly-* in }   L^\infty_{T_0} H^{1.5+\delta}
     ,
     \\&
     a^{(r)}_t \rightharpoonup a_t \text{ weakly-* in }   L^\infty_{T_0} H^{1.5+\delta}
     ,
  \end{split}
  \label{EQ76}
  \end{align}
where we also used \eqref{EQ03}--\eqref{EQ05} and returned to indicating~$r$. Moreover, by the Aubin-Lions lemma it follows that, on the same time interval, 
  \begin{align}
  \begin{split}
     &
     v^{(r)} \to v \text{ strongly in }   C_{T_0} H^{s+1}
     ,
     \\&
     \chi \eta^{(r)} \to \chi \eta \text{ strongly in }   C_{T_0} H^{s+1.5}
     ,
     \\&
     \eta^{(r)} \to \eta \text{ strongly in }   C_{T_0} H^{s+1}
     ,
     \\&
     \chi^2 a^{(r)} \to \chi^2 a \text{ strongly in }   C_{T_0} H^{s+0.5}
     ,
     \\&
     a^{(r)} \to a \text{ strongly in }   C_{T_0} H^{s}
     ,
  \end{split}
  \label{EQ77}
  \end{align}
for any $s < 1.5+\delta$. In order to show that the limit  $(v,q,\eta,a)$ is indeed a solution of the Euler system \eqref{EQ03}--\eqref{EQ07}, we note that we aim to take the limit in the weak formulation of the problem. We note that  the linear terms converge, including the initial condition for $\eta$, as well as the boundary terms by using trace estimates. Recalling the discussion on \eqref{EQ09}--\eqref{EQ11}, it follows that the initial condition for $v$ is satisfied by construction. We thus only need to verify convergence of the nonlinear terms. To this end, we let $\phi \in C_c^{\infty}(\Omega \times (0,T_0))$ and observe that
  \begin{align}
     (a_{ki}^{(r)} \partial_kq^{(r)} - a_{ki} \partial_k q, \phi)
      \lec
       \Vert q^{(r)}\Vert_{L^\infty_{T_0} H^{2.5+\delta}} \Vert a^{(r)} - a\Vert_{L^1_{T_0} {L^1}} 
       \Vert \phi\Vert_{L^\infty_{T_0} L^{\infty}} + (\partial_k q^{(r)} - \partial_k q, a_{ki} \phi),
  \llabel{EQ78}
  \end{align}
which goes to zero, as $r \to 0$, by \eqref{EQ76} and~\eqref{EQ77}.
Next, for the divergence-free condition, we obtain
  \begin{align}
     (a_{ik}^{(r)} \partial_i v_k^{(r)} - a_{ik} \partial_i v_k, \phi) 
      \lec 
       (\Vert a^{(r)}\Vert_{L^\infty_{T_0} H^{1.5+\delta}} \Vert v^{(r)} - v\Vert_{L^\infty_{T_0} H^{1}}
        + \Vert v\Vert_{L^\infty_{T_0} H^{2.5+\delta}} \Vert a^{(r)} - a\Vert_{L^1_{T_0} L^1}
         ) \Vert \phi\Vert_{L^\infty_{T_0} L^{\infty}}, 
  \label{EQ79}
  \end{align}
which also vanishes in the limit. Finally, in the view of \eqref{EQ77} and \eqref{EQ16}, passing to the limit in the right-hand side of \eqref{EQ05} is straightforward, concluding the proof of the existence part of Theorem~\ref{T01}.

It thus remains to verify Steps~3a and~3b. To this end we shall frequently use the Sobolev product rule 
  \begin{align}
   \Vert fg\Vert_{W^{s,p}}
   \lec
   \Vert f\Vert_{W^{s,p_1}}
   \Vert g\Vert_{L^{p_2}}
   +
   \Vert f\Vert_{L^{q_1}}
   \Vert g\Vert_{W^{s,q_2}}
   ,
  \label{EQ16}
  \end{align}
where $s\ge 0$, $p_i, q_i \in [1,\infty]$ and $p \in (1,\infty)$ are such that 
$\frac{1}{p_1} + \frac{1}{p_2} = \frac{1}{p} = \frac{1}{q_1} + \frac{1}{q_2}$. We also recall the commutator estimate~\cite{KP,Li}
  \begin{align}
   \| J(fg) - f Jg \|_{L^2} 
   \lec 
   \| f \|_{W^{s,p_1}} \| g \|_{L^{p_2}} 
   + 
   \| f \|_{W^{1,q_1}} \| g \|_{W^{s-1,q_2}}
   ,
  \label{EQ17}
  \end{align}
for $s\geq 1$ and $\frac{1}{p_1}+\frac1{p_2}= \frac{1}{q_1}+\frac1{q_2}= \frac12$,
and the double commutator estimate
  \begin{align}
   \| J(fg) - f Jg - g Jf \|_{L^p} 
   \lec 
   \| f \|_{W^{1,p_1}} \| g \|_{W^{s-1,p_2}}
   + 
   \| f \|_{W^{s-1,q_1}} \| g \|_{W^{1,q_2}}
   ,
  \label{EQ18}
  \end{align}
for $s\geq 1 $, $\frac{1}{p_1}+\frac1{p_2}=\frac{1}{q_1}+\frac1{q_2}= \frac1p$, $p\in (1,p_1)$, and $p_2,q_1,q_2<\infty$, 
where $J$ is a nonhomogeneous differential operator in $(x_1,x_2)$ of order $s\geq 0$. Finally, we recall the Brezis-Bourgonion \cite{EQ19} inequality
  \begin{align}
   \Vert f\Vert_{{s}} \lec
    \Vert f\Vert_{L^{2}} + \Vert \curl f\Vert_{{s-1}} + \Vert \div f\Vert_{{s-1}} 
     + \Vert \nabla_2 f\Vert_{H^{s-1.5}(\partial \Omega)}, \qquad s \ge 1
     .
  \label{EQ19}
  \end{align}

\subsection{Proof of Step~3a}\label{sec_step3a}

Here we show \eqref{EQ15a} and~\eqref{EQ15b}. 

First, recall the Cauchy invariance
  \begin{align}
     \epsilon_{ijk}\partial_j v_m \partial_k \eta_m = (\omega_0)_i
   \comma i=1,2,3
     ;
  \label{EQ22}
  \end{align}
see \cite[Appendix]{KTV} for a short proof.
From \eqref{EQ22}, we obtain
  \begin{align}
   \begin{split}
    \nabla ((\curl \eta)_i)
    &=
    \epsilon_{ijk} (\delta_{km} - \partial_k \eta_m)\partial_j \nabla (\eta_m) 
    + t  \nabla (\omega_0)_i
    +
   2 \int_0^t
     \epsilon_{ijk} \partial_k v_m \partial_j \nabla \eta_m
        \,\d s
   ,
  \end{split}
   \label{EQ23}
  \end{align}
for $i=1,2,3$.
Applying the fractional product rule \eqref{EQ16} on \eqref{EQ23}, we get
\begin{align}
  \begin{split}
   \Vert \nabla\curl\eta\Vert_{{1+\delta}}
   &\lec
   \epsilon
   \Vert \eta\Vert_{{3+\delta}}
   + t \Vert \nabla \omega_0\Vert_{{1+\delta}} 
   +
   \int_0^t
    \Vert v\Vert_{{2.5+\delta}}
    \Vert \eta\Vert_{{3+\delta}}
   \,\d s
   ,
  \end{split}
   \label{EQ24}
  \end{align}
since
$   \Vert I - \nabla \eta\Vert_{{1.5+\delta}}\lec\epsilon$.
Note that $\omega_0$ actually stands for $\omega_0^{(r)}$ and is finite,
with the bound depending on~$r$.
For the divergence, we use the Fundamental Theorem of Calculus
to obtain
  \begin{align}
  \begin{split}
   \partial_{l} \div \eta
   &=
    (\delta_{kj} - a_{kj}) \partial_{lk} \eta_j
     + \int_0^t
     (\partial_t a_{kj} \partial_{lk} \eta_j
       + \partial_l a_{kj} \partial_k v_j
     ) \,\d s
  ,
  \end{split}
   \label{EQ25}
  \end{align}
for $l=1,2,3$,
where we also used $  a_{kj}\partial_{kl} v_j=-\partial_l a_{kj} \partial_k v_j$
and that $\partial_{lk}\eta$ vanishes at the initial time.
Therefore,
\begin{align}
  \begin{split}
   \Vert \nabla \div  \eta\Vert_{{1+\delta}}
   &\lec
   \Vert I-a\Vert_{{1.5+\delta}}
   \Vert \eta\Vert_{{3+\delta}}
   +
    \int_0^t
      (
        \Vert a_t\Vert_{{1.5+\delta}}
     \Vert \eta\Vert_{{3+\delta}}
        +
        \Vert a\Vert_{{2+\delta}}
        \Vert v\Vert_{{2.5+\delta}}
      )
      \,\d s
  .
  \end{split}
   \label{EQ26}
  \end{align}
Using the curl estimate \eqref{EQ23}, the divergence estimate \eqref{EQ26},
as well as the boundary estimate
  \begin{equation}
   \Vert \eta\Vert_{H^{2.5+\delta}(\Gamma_1)}
   =
   \Vert \chi\eta\Vert_{H^{2.5+\delta}(\Gamma_1)}
   \lec
   \Vert \chi\eta\Vert_{{3+\delta}}
   \llabel{EQ121}
  \end{equation}
in the Brezis-Bourgonion inequality \eqref{EQ19}, we obtain the required control of $\| \eta \|_{3+\delta}$, which proves~\eqref{EQ15a}.

As for \eqref{EQ15b}, we first recall, from \cite[p.~644]{KO} the pressure equation
  \begin{align}
  \begin{split}
    \Delta q
       &=
        \partial_j \bigl( (\delta_{jk} - a_{ji}a_{ki} ) \partial_k q \bigr) + \partial_j (a_{ji} a_{ki}  \partial_k q)
   =
     \partial_j \bigl((\delta_{jk} - a_{ji}a_{ki} ) \partial_k q\bigr)
      + \partial_t a_{ji} \partial_jv_i
    ,
  \end{split}
  \label{EQ27}
  \end{align}
with the boundary conditions
  \begin{equation}
    q=0
    \inon{on $\Gamma_1$}
   \label{EQ28}
  \end{equation}
and
  \begin{equation}
   \partial_3 q = (\delta_{k3} - a_{k3} )\partial_k q
   \inon{on $\Gamma_0$}
   .
   \label{EQ29}
  \end{equation}
Applying the elliptic regularity, using \eqref{EQ15a},
and recalling \eqref{EQ13}--\eqref{EQ14}, we obtain $\| q \|_{3+\delta } \lec \| v \|_{2.5+\delta }^2$, as required.

\subsection{Proof of Step~3b}\label{sec_step3b}
Here we show that if 
\eqnb\label{K}
 \| v_0 \|_4 + \| \omega_0 \|_{3.5} + \sup_{t\in [0,T_0]} \left( \Vert \eta\Vert_{{3+\delta}} +
   \Vert q\Vert_{{3+\delta}}\right) \leq K 
   ,
\eqne
for some $K\geq1$ and \eqref{EQ13} holds for some sufficiently small $\varepsilon >0$, i.e.,
\[ \sup_{t\in [0,T_0]} \left( \| I -a  \|_{1.5+\delta }+\| I -aa^T  \|_{1.5+\delta }+ \| I-\na \eta \|_{1.5+\delta } \right)
  \leq \varepsilon,\]
  then
    \begin{align}
     \Vert \eta (t) \Vert_{{4.5}}^2 + \Vert v (t) \Vert_{{4}}^2
      \lec_K
       C_{\Vert \eta (0) \Vert_{{4.5}}, \Vert v (0) \Vert_{{4}}}
        + \int_0^t  \left( \Vert \eta (s) \Vert_{{4.5}} + \Vert v (s) \Vert_{{4}}+1\right)^2 \,\d s
        ,
  \label{EQ74}
  \end{align}
  where, as indicated, the implicit constant may depend on $K$, proving the claim of Step~3b.
  We emphasize that, although the constant in \eqref{EQ74} may depend on $K$, the choice of $\varepsilon$ may not (as $K$ may grow as $r\to 0^+$, and we need $T_0$ (from Lemma~\ref{L01}) to remain independent of $r$). Thus we need to be careful in the estimates below, so that whenever we need to absorb a term (say, $\varepsilon \| \eta \|_{4.5}$ by the left-hand side) the implicit constant next to the term may not depend on~$K$. We thus keep dependence on $K$ (i.e., the terms $\| \eta \|_{3+\delta }$ and $\|q \|_{3+\delta }$ from Step~3a) explicit, and only use the symbol  ``$\lec$'' to denote an inequality ``$\leq C$'' where $C>0$ may depend on the terms from the a~priori estimate \eqref{EQ15}, i.e., $\| v\|_{2.5+\delta }$, $\| q\|_{2.5+\delta }$, $\| \eta \|_{2.5+\delta }$, and $\| \chi \eta \|_{3+\delta }$ (for all $t\in [0,T_0]$; recall \eqref{EQ15_cons}), which are independent of~$K$.

First, note that, using \eqref{EQ05}, we
have
\begin{equation}
   \Vert a\Vert_{{\sigma}}
   \lec
   \Vert \eta\Vert_{{2.5+\delta}}  \Vert \eta\Vert_{{\sigma+1}}   
   \lec
   \Vert \eta\Vert_{{\sigma+1}} 
   \lec 1
   \comma \sigma\in[0,1.5+\delta]
   \llabel{EQ30}
  \end{equation}
and
  \begin{equation}
   \Vert a\Vert_{{\sigma}}
   \lec
   \Vert \eta\Vert_{{2.5+\delta}}  \Vert \eta\Vert_{{\sigma+1}}   
   \lec
   \Vert \eta\Vert_{{\sigma+1}}
   \comma \sigma\in[1.5+\delta,3.5]
   ,
   \label{EQ31}
  \end{equation}
where we used the multiplicative Sobolev inequality~\eqref{EQ16}. 
Moreover, upon differentiating \eqref{EQ05} in time and using \eqref{EQ04}, it follows that
  \begin{align}
   \Vert a_t\Vert_{{\sigma}}
   \lec
   \Vert v\Vert_{{\sigma+1}} + \Vert \eta\Vert_{{\sigma+1}}
   \lec 1 + \Vert v\Vert_{{\sigma+1}} + \int_0^t \Vert v\Vert_{{\sigma+1}} \,\d s 
   \comma \sigma\geq0
  .
  \label{EQ32}
  \end{align}
Also, we may use the identity $I-a = \int a_t = \int \nabla v \nabla \eta$,
written symbolically,
which  implies that
   \begin{align}
     \begin{split}
      \Vert I-a\Vert_{\sigma}
      &\lec 
        \int_{0}^{t}
   \bigl(
    \Vert v\Vert_{{\sigma+1}}\Vert \eta\Vert_{{1.5+\delta}} 
          +\Vert \eta\Vert_{{\sigma+1}}\Vert v\Vert_{{1.5+\delta}}
   \bigr)\,\d s
    \\&
    \lec
        \int_{0}^{t}
   \bigl(
    \Vert v\Vert_{{\sigma+1}}
          +\Vert \eta\Vert_{{\sigma+1}}
   \bigr)\,\d s
        ,
  \end{split}
        \label{EQ33}
   \end{align}
for $\sigma \geq0$. We show \eqref{EQ74} by first using div-curl estimates to bound $\| \eta \|_{4.5}$ and  $\| v \|_4$ using the  $L^2$ norms of fourth order tangential derivatives of $v$ and $\eta_3$ (Section~\ref{sec_divcurl}). We then derive an estimate for $\| q\|_{4.5}$ in Section~\ref{sec_p}. Using the pressure estimate we then  control the fourth order  tangential derivatives of $v$ in Section~\ref{sec_tan}. We then use this control in  Section~\ref{sec_com}, where we estimate the fourth order tangential derivatives of $\eta$ and complete the proof of \eqref{EQ74}.

\subsubsection{Div-Curl Estimates}\label{sec_divcurl}
We first aim to show that
  \begin{align}
    \begin{split}
   \Vert \eta\Vert_{{4.5}}
   &\lec
   K + K \int_{0}^{t} (\Vert \eta\Vert_{{4.5}}  +   \Vert v\Vert_{{4}} )\,\d s
   +
     \sum_{|\alpha| = 4; \alpha_3=0}
     \Vert \partial^{\alpha} \eta_3\Vert_{L^{2}(\Gamma_1)} 
  \end{split}
  \label{EQ34}
  \end{align}
and 
  \begin{align}
   \Vert v\Vert_{{4}} \lec
     K
     + \int_0^t
\Vert v\Vert_{{4}}
\,\d s 
      + \sum_{|\alpha| = 4; \alpha_3=0} \Vert \partial^{\alpha} v\Vert_{L^{2}}
  .
  \label{EQ35}
  \end{align}
  
In order to obtain \eqref{EQ34}, we first apply the fractional Sobolev inequalities \eqref{EQ16} on the formula \eqref{EQ23} for $\na \curl \eta$ and write
  \begin{align}
  \begin{split}
   \Vert \nabla\curl \eta\Vert_{{2.5}} 
   &\lec
   \Vert I - \nabla \eta\Vert_{{1.5+\delta}}
   \Vert  \eta\Vert_{{4.5}}
   +
   \Vert I - \nabla \eta\Vert_{{3}}
   \Vert \eta\Vert_{{3}}
   \\&\indeq
   + t \Vert \nabla \omega_0\Vert_{{2.5}} 
   +
   \int_0^t
   (
    \Vert \nabla v\Vert_{{1.5+\delta}}
    \Vert \eta\Vert_{{4.5}}
    +
    \Vert  v\Vert_{{4}}
    \Vert \eta\Vert_{{3}}
   )   
   \,\d s
   \\
   &\lec
   \varepsilon \| \eta \|_{4.5}
   + K (\| \eta \|_4+1)
   +
   K
   \int_0^t
   (
    \Vert \eta\Vert_{{4.5}}
    +
    \Vert  v\Vert_{{4}}
   )   
   \,\d s
  \\
   &\lec
   \varepsilon \| \eta \|_{4.5}
   + C_\varepsilon K
   +
   K
   \int_0^t
   (
    \Vert \eta\Vert_{{4.5}}
    +
    \Vert  v\Vert_{{4}}
   )   
   \,\d s,
  \end{split}
  \label{EQ38}
  \end{align}
where we used the interpolation
  $\| \eta \|_{4}\lec \| \eta \|_{4.5}^{3/4} \| \eta \|_{2.5}^{1/4} \leq \varepsilon \| \eta \|_{4.5} + C_\varepsilon \| \eta \|_{2.5}\lec \varepsilon\Vert \eta\Vert_{4.5}+C_{\varepsilon}$ and \eqref{K} in the last inequality.
  
For $\div \eta $ we estimate  \eqref{EQ25} in $H^{2.5}$ to obtain
  \begin{align}
  \begin{split}
   \Vert \nabla \div \eta\Vert_{{2.5}}
   &\lec
   \Vert I-a\Vert_{{1.5+\delta}}
   \Vert \eta\Vert_{{4.5}}
    +
   \Vert I-a\Vert_{{3}}
   \Vert \eta\Vert_{{3}}
   \\&\indeq
    +
    \int_0^t
      \Bigl(
        \Vert a_t\Vert_{{1.5+\delta}}
        \Vert \eta\Vert_{{4.5}}
        +
        \Vert a_t\Vert_{{3}}
        \Vert \eta\Vert_{{3}}
       + \Vert a\Vert_{{2}}
         \Vert v\Vert_{{4}}
       + \Vert a\Vert_{{3.5}}
        \Vert v\Vert_{{2.5+\delta}}
      \Bigr)
      \,\d s\\
      &\lec \varepsilon \| \eta \|_{4.5}
    + C_{\varepsilon} K  
    +
   K \int_0^t
      \Bigl(
        \Vert \eta\Vert_{{4.5}}
       + 
         \Vert v\Vert_{{4}}
      \Bigr)
      \,\d s
   ,
  \end{split}
  \label{EQ39}
  \end{align}
where we used \eqref{EQ31}, \eqref{EQ32}, and \eqref{EQ33} in the last line;
we have also treated 
the term $ \Vert I-a\Vert_{{3}}    \Vert \eta\Vert_{{3}}$
analogously to
$
   \Vert I - \nabla \eta\Vert_{{3}}
   \Vert \eta\Vert_{{3}}
$ in~\eqref{EQ38}.
Applying these in the Brezis-Bourgonion inequality \eqref{EQ19} gives 
   \begin{align}
  \begin{split}
   \Vert  \eta\Vert_{{4.5}}
     &\lec
   1
   +
   \varepsilon \Vert \eta\Vert_{{4.5}}
   +C_\varepsilon K
   + K\int_{0}^{t} (\Vert \eta\Vert_{{4.5}}  +   \Vert v\Vert_{{4}} )\,\d s
   + \sum_{|\alpha| = 4; \alpha_3=0} \Vert \partial^{\alpha} \eta_3\Vert_{L^{2}(\Gamma_1)} 
  .
  \end{split}
  \llabel{EQ41}
  \end{align}
We may now absorb the term with $\varepsilon$ by the left-hand side to obtain~\eqref{EQ34}.

To prove \eqref{EQ35}, we first observe that the Cauchy invariance \eqref{EQ22} gives
  \begin{align}
  (\curl v)_i = \epsilon_{ijk} (\delta_{km} - \partial_k \eta_m) \partial_j v_m + (\omega_0)_i
   \comma i=1,2,3
   ,
  \label{EQ42}
  \end{align}
and the divergence-free condition \eqref{EQ03} implies
  \begin{align}
  \div v = (\delta_{ji} - a_{ji}) \partial_j v_i.
  \label{EQ43}
  \end{align}
Using the Brezis-Bourgonion inequality \eqref{EQ19} again, we obtain
\begin{align}
  \Vert v\Vert_{{4}}
  \lec
  \Vert v\Vert_{L^{2}}
  + \sum_{i}\Vert \epsilon_{ijk} (\delta_{km} - \partial_k \eta_m) \partial_j v_m\Vert_{{3}}
  + 
  \Vert (\delta_{ji} - a_{ji}) \partial_j v_i\Vert_{{3}} + \Vert \nabla_2 v_3\Vert_{H^{2.5}(\Gamma_1)}
  + \Vert \omega_0\Vert_{{3}}
  .   
  \label{EQ44}
  \end{align}
For the boundary term, we use the trace estimate to write
  \begin{align}
  \Vert \nabla_2 v_3\Vert_{H^{2.5}(\Gamma_1)} \lec 
  \Vert v\Vert_{L^{2}}
    + \sum_{|\alpha| = 3; \alpha_3=0}
      \Vert \partial^\alpha \nabla v_3\Vert_{L^{2}}
  .
  \llabel{EQ45}
  \end{align}
Together with
  \begin{align}
  \partial_3 v_3 = (\delta_{ji} - a_{ji}) \partial_j v_i - \partial_1 v_1 - \partial_2 v_2,
  \label{EQ46}
  \end{align}
it follows by Lemma~\ref{L01} that
if
$|\alpha|=3$ and
$\alpha_3=0$, then
\begin{align}
     \Vert \partial^\alpha \nabla v_3\Vert_{L^{2}}
     \lec 
      \Vert (I - a) \nabla v\Vert_{{3}}
       + \sum_{|\alpha'| = 4;\alpha'_3=0}
       \Vert \partial^{\alpha'}  v\Vert_{L^{2}} 
   .
   \llabel{EQ02}
  \end{align}
Substituting this in \eqref{EQ44} gives
  \[
  \begin{split}
   \Vert v\Vert_{{4}}
  &\lec
  K
  + \| (I-\na \eta )\na v \|_3+ \| (I-a)\na v \|_3 +\sum_{|\alpha'| = 4;\alpha'_3=0}
       \Vert \partial^{\alpha'}  v\Vert_{L^{2}} \\
       &\lec K+ \varepsilon \| v \|_4 + (\| I-a \|_3+\| I-\na \eta \|_3 ) \|v\|_{2.5+\delta } +\sum_{|\alpha'| = 4;\alpha'_3=0}
       \Vert \partial^{\alpha'}  v\Vert_{L^{2}} .
\end{split}
  \]
Absorbing $\varepsilon \| v \|_4 $, recalling the definition \eqref{EQ04} of $\eta$, the fact $\| v \|_{2.5+\delta }\lec 1$, along with \eqref{EQ31} gives \eqref{EQ35}, as required.

\subsubsection{Pressure Estimates}\label{sec_p}
Here we show that
  \begin{align}
     \Vert q\Vert_{{4.5}} \lec
      1+ \Vert \eta\Vert_{{4.5}}
      .
  \label{EQ48}
  \end{align}
Using the $H^{4.5}$ elliptic regularity \cite[Chapter~2]{LM} on the Poisson problem \eqref{EQ27}--\eqref{EQ29} for $q$ gives 
  \begin{align}
  \begin{split}
     \Vert q\Vert_{{4.5}}
     &\lec
      (\Vert I - aa^T\Vert_{{1.5+\delta}}
       + \Vert I - a\Vert_{{1.5+\delta}}
      )
      \Vert q\Vert_{{4.5}} 
       +
       (\Vert I - aa^T\Vert_{{3.5}}
       + \Vert I - a\Vert_{{3.5}}
       )
       \Vert  q\Vert_{{2.5+\delta}}
    \\&
    \lec
      \varepsilon
      \Vert q\Vert_{{4.5}} 
       +
       1 + \Vert \eta\Vert_{{4.5}}
   ,
  \end{split}
   \llabel{EQ49}
  \end{align}
  where we used \eqref{EQ13}, \eqref{EQ31} and recalled from \eqref{EQ15_cons} that $\Vert q\Vert_{{2.5+\delta}}\lec 1$
in the last step. Absorbing $\varepsilon\Vert q\Vert_{{4.5}}$ by the left side, we obtain~\eqref{EQ48}.

\subsubsection{Tangential Estimates}\label{sec_tan}
In this section, we prove that
\begin{align}
   \Vert D\partial^\alpha  v\Vert_{L^{2}}^2 + \Vert \tau a_{3i} D \partial^\alpha \eta_i \Vert_{L^{2}(\Gamma_1)}^2 
   \lec_K 1
   + \int_0^t  Q(s)\,\d s
   \comma |\alpha|=3
   \commaone \alpha_3=0
  ,   
  \label{EQ50}
  \end{align}
where we set
\[
Q=(\Vert \eta\Vert_{{4.5}}  +   \Vert v\Vert_{{4}} + 1)^2
\]
for brevity, while the symbols $D$ and $\tau$ denote
  \begin{align}
  \begin{split}
     Df(x) &= \frac{1}{h}(f(x+he_l)-f(x))
   \comma l=1,2
  \end{split}
   \llabel{EQ37}
  \end{align}
and
  \begin{align}
  \begin{split}
     \tau f(x)&= f(x+he_l)
   \comma l=1,2
   ,
  \end{split}
  \llabel{EQ51}
  \end{align}
for
$h \in (0,1]$;
observe that we only consider tangential directions when using $D$ or~$\tau$.
Recall that the product rule for $D$ reads
  \begin{align}
  D(fg) = Df \tau g + fDg
  .
  \label{EQ52}
  \end{align}
We note that we use the difference quotient $D$, since we only have control of $\eta \in H^{4.5}$ (recall that we are showing \eqref{EQ74}), while some of the integrals resulting from estimates of $\| \nabla_2 \p^\alpha v \|_{L^2}$  would involve five derivatives of~$\eta$.
We thus replace tangential gradient $\nabla_2$ by the tangential difference quotient $D$, so that our manipulations on the integrals are justified, and we take $h\to 0^+$ at the end. We also point out that we the inequality
\[
\| D f \|_s \lec \| \na f \|_s,\qquad s\geq 0
,
\]  
which is a consequence of the Fundamental Theorem of Calculus,
in numerous estimates.
  
Now, fix $\alpha$ such that $|\alpha|=3$ and $\alpha_3=0$.
We apply $D\p^\alpha$ to \eqref{EQ03}, test it with 
$D\partial^\alpha  v_i$ and use the product rule \eqref{EQ52} to obtain
\begin{align}
  \begin{split}
     \frac{1}{2}\frac{\d}{\d t}\Vert D\partial^\alpha v\Vert_{L^{2}}^2
     &=
     - 
     \int  
      (\partial^\alpha (a D \nabla q) - a \partial^\alpha D\nabla q) 
       D\partial^\alpha v
     -
     \int  
      a_{ki}\partial^\alpha D\partial_k q
       D\partial^\alpha v_i
       \\&\indeq
     -
     \int  
      (\partial^\alpha (Da \tau \nabla q) - \partial^\alpha D a\, \tau  \nabla q)
       D\partial^\alpha v
     -
     \int  
      \partial^\alpha D a_{ki}\tau \partial_k q
       D\partial^\alpha v_i
     \\&
     = I_1 + I_2 + I_3 + I_4
  \end{split}  
  \label{EQ53}
  \end{align}
(using a symbolic notation where the indices are not important),
as $D$ and $\partial^\alpha$ commute.
We estimate the first and third terms using \eqref{EQ17}, obtaining
  \begin{align}
  \begin{split}
  &
     I_1 \lec 
     \Vert \partial^\alpha (a D  \nabla q) - a \partial^\alpha D \nabla q\Vert_{L^{2}}
     \Vert D \partial^\alpha v\Vert_{L^{2}}
     \lec 
     (\Vert  a\Vert_{{3.5}} 
     + 
     \Vert q\Vert_{{4.5}}) \Vert v\Vert_{{4}}
     \lec_K Q
  \end{split}
  \label{EQ54}
  \end{align}
and
  \begin{align}
  \begin{split}
  &
     I_3 \lec 
     \Vert \partial^\alpha (Da \tau \nabla q) - \partial^\alpha D a \tau \nabla q\Vert_{L^{2}}
     \Vert D \partial^\alpha v\Vert_{L^{2}}
     \lec 
     (\Vert  q\Vert_{{4.5}} 
     + \Vert  a\Vert_{{3.5}}
     )
     \Vert v\Vert_{{4}}\lec_K Q
  ,
  \end{split}
  \label{EQ55}
  \end{align}
where we used  \eqref{EQ31} and \eqref{EQ48} in the last
inequalities, respectively.
For $I_2$, we integrate by parts in $x_k$ the term involving $\partial_k q$ which yields
  \begin{align}
     I_2 = 
      \int  
      a_{ki} \partial^\alpha D q D \partial^\alpha \partial_k v_i
  \llabel{EQ56}
  ,
  \end{align}
since the Piola Identity $\partial_k a_{ki} = 0$, for $i=1,2,3$, eliminates the term involving the partial derivative of $a_{ki}$, while \eqref{EQ07} and $a_{3i}|_{\Gamma_0} = 0$, for $i=1,2$ removes the boundary terms
on $\Gamma_1$ and $\Gamma_0$, respectively.
Using
\eqref{EQ03}$_2$, we may rewrite $I_2$ using the product rule \eqref{EQ52} to obtain
\begin{align}
  \begin{split}
      I_2
      &=
      -
      \int  
      \partial^\alpha( D a_{ki} \tau \partial_k v_i)
      D\partial^\alpha q
      -
      \int  
      \Bigl(
        \partial^\alpha (a_{ki} D \partial_k v_i) - a_{ki} \partial^\alpha D \partial_k v_i
      \Bigr) 
      D\partial^\alpha q
     \\&
      = I_{21} + I_{22}
   .
  \end{split}
  \llabel{EQ57}
  \end{align}
We estimate $I_{22}$ similarly to \eqref{EQ54} and \eqref{EQ55} using~\eqref{EQ17}, leading to
  \begin{align}
  \begin{split}
     I_{22}
     &\lec
     \Vert \partial^\alpha (a D \nabla v) - a \partial^\alpha D \nabla v\Vert_{L^{3/2}}
      \Vert D\partial^\alpha q\Vert_{L^{3}}
     \lec
     (\Vert a\Vert_{{3}} \Vert D\nabla v\Vert_{L^{6}} 
      + \Vert  a\Vert_{W^{1,6}} \Vert D\nabla v\Vert_{{2}}) \Vert q\Vert_{{4.5}}
      \\&
      \lec 
     ( \Vert  a\Vert_{{3.5}}
     +
     \Vert v \Vert_{{4}} )
     \Vert q\Vert_{{4.5}} \lec_K Q
  .
  \end{split}
   \llabel{EQ59}
  \end{align}
On the other hand, for $I_{21}$ we have
  \begin{align}
     I_{21} =
      -
    \int
      \Lambda^{-1/2}
      \partial^\alpha( D a_{ki} \tau \partial_k v_i)
      \Lambda^{1/2}
      D\partial^\alpha q
      \lec_K
      (
      \Vert a\Vert_{{3.5}}  
      +
       \Vert v\Vert_{{4}}
       )
       \Vert q\Vert_{{4.5}}
     \lec_K Q
       ,
  \llabel{EQ60}
  \end{align}
where $\Lambda \coloneqq (I_2-\Delta_2)^{\frac{1}{2}}$, using again a symbolic notation, i.e.,
omitting indication of indices.
In the second step, we used \eqref{EQ15b} and
\eqref{EQ31}, \eqref{EQ48} in the last.

Thus, in order to prove \eqref{EQ50}, it remains to estimate~$I_4$. To do so, we first apply $\partial^\alpha$ to the equality $a \nabla \eta = I$ and solve for $\partial^\alpha a$, thereby obtaining 
$\partial^{\alpha} a_{ki} = - a_{kj} \partial^{\alpha} \partial_{l} \eta_j a_{li}
 - Ra$, where
\begin{align}
  R = \partial^\alpha (a \nabla \eta) - \partial^\alpha a \nabla \eta - a \partial^\alpha \nabla \eta
  .
  \label{EQ61}
\end{align}
Next, we apply $D$ and use \eqref{EQ52}, obtaining
\begin{align} 
  D\partial^\alpha a_{ki} 
   =
    -\tau(a_{kj}a_{li})D\partial^\alpha \partial_l \eta_j
     - D(aa)\partial^\alpha \nabla \eta
      - D(Ra)
      ,
      \llabel{EQ62}
\end{align}
from where, rearranging the terms, 
  \begin{align}
  \begin{split}
      I_4
      &=
      \int  
       \bigl(D(aa) \partial^\alpha \nabla \eta + D(Ra)\bigr)
         \tau \nabla q D \partial^\alpha v 
        \\&\indeq
        + 
      \int
       \tau (a_{kj}a_{li}) D \partial^\alpha \partial_l \eta_j
        \tau \partial_k q D \partial^\alpha v_i
      =: I_{41} + I_{42}.
  \end{split} 
  \label{EQ63}
  \end{align}
Observe that for $R$, defined in \eqref{EQ61}, we have 
\begin{align}
  \begin{split}
   \Vert R\Vert_{{1}}
   &\lec
     \Vert a\Vert_{{3.5}}\Vert  \eta\Vert_{{3+\delta}} 
      + \Vert a\Vert_{{2}}\Vert \eta\Vert_{{4.5}}
    \lec
    K
    (
     \Vert a \Vert_{{3.5}}
     + 
     \Vert \eta \Vert_{{4.5}}
    )
    .
  \end{split}
  \label{EQ64}
\end{align}
Using \eqref{EQ64}, the first term in \eqref{EQ63} is estimated by 
\begin{align}
  I_{41}
   \lec 
    (\Vert D(aa)\Vert_{L^{6}}\Vert \partial^\alpha \nabla \eta\Vert_{L^{3}} + \Vert D(Ra)\Vert_{L^{2}}) 
     \Vert D\partial^\alpha v\Vert_{L^{2}}
      \lec_{K}
       (\Vert \eta\Vert_{{4.5}} + \Vert R\Vert_{{1}}) \Vert v\Vert_{{4}} \lec_K Q
       ,
      \llabel{EQ65}
\end{align} 
as required.

Next, for $I_{42}$, we use integration by parts in $l$ which yields
  \begin{align}
  \begin{split}
     I_{42}
     &= -
     \int  
      \Bigl(
      \partial_l \tau (a_{kj}a_{li}) \partial^\alpha D \eta_j 
      \tau \partial_k q
      +
      \tau (a_{kj}a_{li}) \partial^\alpha D \eta_j 
      \tau\partial_{lk} q
      \Bigr) D \partial^\alpha  v_i
      \\&\indeq
     + 
     \int_{\Gamma_0\cup\Gamma_1}
     \tau (a_{kj}a_{li}) \partial^\alpha D \eta_j 
     \tau \partial_k q D \partial^\alpha  v_i N_{l}
     \,\d\sigma
     \\&\indeq
     -
     \int  
      \tau (a_{kj}a_{li}) \partial^\alpha D \eta_j 
      \tau \partial_k q D \partial^\alpha \partial_l v_i
     = I_{421} + I_{422} + I_{423}.
  \end{split} 
  \llabel{EQ66}
  \end{align}
The first term $I_{421}$ consists of lower order terms that may be
estimated as
  \begin{align}
  \begin{split}
   I_{421}
   &\lec
   \| \na \tau(a^2) \na^3D \eta \na q\|_{L^2} \| Dv \|_3
   +
   \Vert
   \tau (a^2)
   \nabla^{3}D \eta
   \nabla^2 q
   \Vert_{L^2}  \| D v \|_3
   \\&
   \lec
   \| a \|_{1.5+\delta }
   \| a \|_{2+\delta }
   \| \eta \|_{4.5}
   \| q \|_{3 }
   \| v \|_4
   \lec_K
   \Vert \eta\Vert_{4.5}    \| v \|_4 \lec Q
  ,
  \end{split}
   \llabel{EQ21}
  \end{align}
where we used~\eqref{EQ31}.
For $I_{423}$, we have
  \begin{align}
  \begin{split}
   &
    \tau (a_{kj}a_{li}) 
    D \partial^\alpha \partial_l v_i
   \\&\indeq
    =
    \partial^\alpha ( \tau(a_{kj}a_{li}) D\partial_{l} v_i)
    -     \bigl(
            \partial^{\alpha}( \tau(a_{kj}a_{li})  D\partial_l v_i)
               - \tau(a_{kj}a_{li}) D\partial^{\alpha} \partial_{l}  v_i
         \bigr)
   \\&\indeq
    =
    \partial^\alpha D( a_{kj}a_{li}  \partial_{l}v_i )
    -
    \partial^\alpha ( D(a_{kj}a_{li})  \partial_{l}v_i )
    -     \bigl(
            \partial^{\alpha}( \tau(a_{kj}a_{li})  D\partial_l v_i)
               - \tau(a_{kj}a_{li}) D\partial^{\alpha} \partial_{l} v_i
          \bigr)
    ,
   \end{split}
   \label{EQ68}
   \end{align}
and we note that the first term on the far right vanishes by the divergence-free condition.
Using \eqref{EQ68}, we get
  \begin{align}
  \begin{split}
    I_{423}
      &=
     \int  
      \partial^\alpha D \eta_j 
      \tau \partial_k q
     \partial^\alpha ( D(a_{kj}a_{li}) \partial_l v_i)
      \\&\indeq
           +
     \int  
      \partial^\alpha D \eta_j 
      \tau (\partial_k q)
     \bigl(
            \partial^{\alpha}( \tau(a_{kj}a_{li})  D\partial_l v_i)
               - \tau(a_{kj}a_{li}) D\partial^{\alpha} \partial_{l} v_i
     \bigr)
     \\&
     \lec
     (\Vert  a\Vert_{{3.5}}  
     +
     \Vert v\Vert_{{4}})
     \Vert \eta\Vert_{{4.5}}
     \lec_K Q,
   \end{split}
   \llabel{EQ69}
  \end{align}
due to~\eqref{EQ31}.

It remains to treat the boundary integral~$I_{422}$.
The integral over $\Gamma_0$ vanishes since
$N = (0,0,-1)$ and $a_{3i} = 0$ on~$\Gamma_0$. On the other hand, note
that  $\partial_k q = 0$ for $k=1,2$ on $\Gamma_1$,
which gives
$I_{422}=     \int_{\Gamma_1}
     \tau (a_{3j}a_{3i}) \partial^\alpha D \eta_j 
     \tau \partial_3 q D \partial^\alpha v_i
     \,\d\sigma
$.
Moreover, we may write
$\tau a_{3i} D \partial^\alpha v_i = 
\partial_t(\tau a_{3i} D \partial^\alpha \eta_i) - \partial_t(\tau a_{3i}) D \partial^\alpha \eta_i$, so that $I_{422}$ becomes
  \begin{align}
  \begin{split}
     I_{422} &=
     \frac{1}{2} \frac{\d}{\d t}
      \int_{\Gamma_1}  
       (\tau a_{3i} \partial^\alpha  D \eta_i
       )^2 \tau \partial_3 q  
        \,\d\sigma
     \\&\indeq
     - \frac{1}{2}
      \int_{\Gamma_1}  
       \bigl((\tau a_{3i} \partial^\alpha  D \eta_i)^2 \tau \partial_3 \partial_t q  
     +
     2
       \tau a_{3j} \partial^\alpha  D \eta_j \tau \partial_3 q \partial_t \tau a_{3i} D \partial^\alpha \eta_i
       \bigr)  
        \,\d\sigma 
     ,
  \end{split}      
  \label{EQ70}
  \end{align}
with the second term in \eqref{EQ70} bounded by
$\Vert \psi q_t\Vert_{{2.5+\delta}}  \Vert \chi \eta\Vert_{{4.5}}^2
   +    \Vert \psi q\Vert_{{3+\delta}}
        \Vert \chi \eta\Vert_{{4.5}}^2
\lec \| \eta \|_{4.5}^2 \lec Q$, due to~\eqref{psiqt},
where we also recalled \eqref{EQ32}, which gives $\Vert a_t\Vert_{L^{\infty}}\lec1$. 
We denote by $I(t)$ the time integral on $(0,t)$ of the first term in~\eqref{EQ70}, and we obtain
\begin{align}
  \begin{split}
     I(t) &= \frac{1}{2}
      \int_{\Gamma_1}  
      (\tau a_{3i} \partial^\alpha  D \eta_i)^2 \tau \partial_3 q  
      \,\d\sigma |_t
       - \frac{1}{2}
        \int_{\Gamma_1}  
         (\partial^\alpha  D \eta_3)^2 \tau \partial_3 q  
        \,\d\sigma |_0
     \le
      -\frac{b}{4} \Vert \tau a_{3i} \partial^\alpha  D \eta_i\Vert_{L^{2}(\Gamma_1)}^2
       + C 
       ,  
  \end{split}     
  \label{EQ71}
  \end{align}
where we used $\tau a_{ki}(0) = \tau \delta_{ki}= \delta_{ki}$ in the second term
and Lemma~\ref{L01} and \eqref{K} in the final step.
We thus obtain \eqref{EQ50} by integrating \eqref{EQ53} in time, combining the estimates for $I_1$--$I_4$ above, and by moving $-\frac{b}{4} \Vert \tau a_{3i} \partial^\alpha  D \eta_i\Vert_{L^{2}(\Gamma_1)}^2$ to the left-hand side.\colb  

\subsubsection{Conclusion of the proof of \eqref{EQ74}}\label{sec_com}
Thanks to the control of the tangential derivatives of $v$ from \eqref{EQ50},
we can now estimate the tangential derivatives of $\eta$ and conclude the proof of~\eqref{EQ74}.

First,  for any multiindex $\alpha$ with $|\alpha|=3$ and $\alpha_3=0$, and for any $h\in(0,1]$ and $l\in\{1,2\}$, we have
  \begin{align}
  \begin{split}
     \Vert D\partial^\alpha \eta_3\Vert_{L^{2}(\Gamma_1)}^2
     &\lec 
     \Vert \tau a_{3i} D \partial^\alpha \eta_i\Vert_{L^{2}(\Gamma_1)}^2
     + \Vert (\delta_{3i} - \tau a_{3i}) D \partial^\alpha \eta_i\Vert_{L^{2}(\Gamma_1)}^2
     \\&
     \lec 
     \varepsilon \Vert  \eta\Vert_{{4.5}}^2 + C_K
     + C_K \int_0^t Q(s)  \d s
  ,
  \end{split} 
  \label{EQ73}
  \end{align}
where we used the tangential estimate \eqref{EQ50} and $I - \tau a = \tau (I - a)$, together with Lemma~\ref{L01} in the last inequality. Since the right-hand sides of \eqref{EQ50} and \eqref{EQ73} are independent of $h$, we may vary $l\in \{ 1,2\}$ and $\alpha$, and take the limit $h\to \infty$ to obtain
\[
  \sum_{|\alpha |=4 ; \alpha_3=0}   \left( \Vert \partial^\alpha \eta_3\Vert_{L^{2}(\Gamma_1)}^2    
     + \| \partial^\alpha v \|_{L^2}\right) \lec 
     \varepsilon \Vert  \eta\Vert_{{4.5}}^2 + C_K
     + C_K \int_0^t Q(s)  \d s
  ,
\]
recalling \eqref{psiqt} that $t\leq T_0 \leq 1$.
Finally, squaring the div-curl estimates \eqref{EQ34} and \eqref{EQ35} and adding the resulting inequalities, we obtain
  \begin{align}
  \begin{split}
     \Vert \eta\Vert_{{4.5}}^2 + \Vert v\Vert_{{4}}^2
      &\lec \varepsilon \|\eta \|_{4.5}^2+
       C_K
        + C_K \int_0^t Q(s)  \, \d s,
  \end{split}
  \llabel{EQ72}
  \end{align}
which gives \eqref{EQ74} after absorbing $\varepsilon \| \eta \|_{4.5}^2 $ by the left-hand side.

\section{Proof of Theorem~\ref{T01}: uniqueness}\label{sec_uniqueness}
Having established the existence in the previous section, we now turn to uniqueness. To this end, let $(v,q,\eta,a)$ and $(\tilde{v},\tilde{q},\tilde{\eta},\tilde{a})$ be two solutions as in \eqref{EQ08} emanating from the same initial datum with at least one of them satisfying~\eqref{EQ20},
and denote by $(V,Q,E,A)$ their difference.
We claim that
\eqnb\label{EQ120}
  \Vert V(t) \Vert_{{1.5+\delta}}^2 + \Vert \chi E(t) \Vert_{{2+\delta}}^2 \lec \int_0^t \left( \| V(s)  \|_{1.5+\delta }^2 + \| \chi E (s) \|_{2+\delta }^2 \right) \d s 
,
\eqne
for all $t\in [0,T]$, where the implicit constant may depend on $T$. Uniqueness then follows by a simple Gronwall argument.

In what follows we allow all implicit constants to depend on $T$ and the norms \eqref{EQ08} of the two solutions at any $t\in [0,T]$. For simplicity we  set
  \begin{equation}
   Y \coloneqq \| V \|_{1.5+\delta } + \| \chi E \|_{2+\delta }
   ,
   \llabel{EQ40}
  \end{equation}
for each given $t$.
In order to obtain \eqref{EQ120}, we shall first use div-curl estimates to obtain
\begin{align}
  \Vert \chi E\Vert_{{2+\delta}} \lec 
  \Vert S E_3\Vert_{L^{2}(\Gamma_1)}
  + \int_0^t Y
  \label{EQ90}
\end{align} 
and
\begin{align}
  \Vert V\Vert_{{1.5+\delta}} 
  \lec
  \Vert V\Vert_{L^{2}}
  + \Vert S (\psi V)\Vert_{L^{2}} 
  + \int_0^t \Vert V\Vert_{{1.5+\delta}}\,\d s  
  ,
  \label{EQ96}
\end{align}
where we set
  \begin{equation}
   S\coloneqq \Lambda^{1.5+\delta } = (I_2 - \Delta_2 )^{1/2}
   ;
   \label{EQ67}
  \end{equation}
see Section~\ref{sec_uni_divcurl} for details.
We then show in Section~\ref{sec_uni_tan} that
\begin{align}
  \Vert S(\psi V)\Vert_{L^{2}}^2 + \Vert SE_3 \Vert^2_{L^{2}(\Gamma_1)}
  \lec \varepsilon \| \chi E \|_{2+\delta}+
  \int_0^t Y^2
  \label{EQ118}
  .
\end{align}
The above three inequalities then imply that
\[
 \Vert V(t) \Vert_{{1.5+\delta}}^2 + \Vert \chi E(t) \Vert_{{2+\delta}}^2 \lec \| V \|_{L^2}^2+ \int_0^t Y^2 ,
\]
which together with the inequality
\eqnb\label{EQ121}
\| V \|_{L^2} \lec \int_0^t Y,
\eqne
which we show in Section~\ref{sec_uni_l2}, gives \eqref{EQ120}, as required.

\subsection{Div-curl estimates for $V$ and $\chi E$}\label{sec_uni_divcurl}

Here we prove \eqref{EQ90} and~\eqref{EQ96}. 

To this end we first recall that  the definition  \eqref{EQ04} of $\eta$ gives that  
\begin{align}
  \Vert E\Vert_{{1.5+\delta}} \lec \int_0^t \Vert V\Vert_{{1.5+\delta}}\,\d s
  ,
  \label{EQ85}
\end{align}
and similarly \eqref{EQ05} implies
\begin{align}
  \Vert \chi^2 A\Vert_{{\sigma}} \lec \Vert \chi E\Vert_{{\sigma +1}}
   \comma   \sigma \in [0,1+\delta ]
   ,
  \label{EQ86}
\end{align}
with the same inequality
when $\chi$ is replaced by~$1$.
In addition, we may differentiate in time
the formula for 
$A$ resulting to \eqref{EQ05} to obtain
\begin{align}
  \Vert A_t\Vert_{{0.5+\delta}} \lec \Vert V\Vert_{{1.5+\delta}} + \Vert E\Vert_{{1.5+\delta}}
  .
  \label{EQ87}  
\end{align}

In order to verify \eqref{EQ90} the divergence-free condition in \eqref{EQ03} and the fact that $\eta_t =v$ allow us write, as in \cite[p.~642]{KO},
\begin{align}
\begin{split}
\div\, \p_l (\chi \eta )
& =  ( \delta_{kj} - a_{kj} ) \p_l \p_k (\chi \eta_j ) 
  \\&\indeq
    + \int_0^t  \bigl( \p_t a_{kj} \p_l \p_k (\chi \eta_j )
    +  a_{kj} ( \p_l \p_k \chi  v_j + \p_l \chi \p_k v_j + \p_k \chi \p_l v_j)  - \chi \p_l a_{kj} \p_k v_j  \bigr) \d s
\\&\indeq
 + \delta_{l3}(\partial_{3}\partial_{3}\chi x_2 + 2\p_3\chi) 
,
\end{split}
   \llabel{EQ80}
\end{align}
from where, taking a difference,
\begin{align}
  \begin{split}
    \div \partial_l (\chi E)
    &=
    -A_{kj} \partial_l \partial_k (\chi \eta_j) 
    + (\delta_{kj} - \tilde{a}_{kj}) \partial_l \partial_k (\chi E_j)
    \\&\indeq
    + \int_0^t
      \Bigl(\partial_t A_{kj} \partial_l \partial_k (\chi \eta_j) 
             + \partial_t\tilde{a}_{kj} \partial_l \partial_k (\chi E_j)
             + A_{kj} (\partial_l \partial_k \chi v_j + \partial_l \chi \partial_k v_j + \partial_k \chi \partial_l v_j)
    \\&\indeq\indeq\indeq\indeq\indeq
             + \tilde{a}_{kj} (\partial_l \partial_k \chi V_j + \partial_l \chi \partial_k V_j + \partial_k \chi \partial_l V_j) 
             - \chi \partial_l A_{kj} \partial_k v_j - \chi \partial_l \tilde{a}_{kj} \partial_k V_j
     \Bigr)\,\d s
    .
  \end{split}
  \llabel{EQ81}
\end{align}
Therefore, similarly to \eqref{EQ39}, we obtain
\begin{align}
  \begin{split}
    \Vert  \nabla \div (\chi E)\Vert_{{\delta}}
    &\lec
    \Vert  A\Vert_{{0.5+\delta}} + \Vert  I-a\Vert_{{1.5+\delta}} \Vert  \chi E\Vert_{{2+\delta}}
    \\& \indeq
    + \int_0^t (
      \Vert  A_t\Vert_{{0.5+\delta}} + \Vert \chi E\Vert_{{2+\delta}} 
   + \Vert A\Vert_{\delta}
    + \Vert \chi \nabla A\Vert_{{\delta}} + \Vert V\Vert_{{1.5+\delta}})\,\d s\\
    &\lec  \varepsilon \| \chi E \|_{2+\delta } + \int_0^t Y
    ,
  \end{split}
  \label{EQ82}
\end{align}
where we also used \eqref{EQ85} to simplify double integrals in time into a single one.
As for the curl, we, as in \cite[p.~642]{KO}, use the Cauchy
invariance to obtain
\begin{align}
  \begin{split}
    \nabla ((\mathrm{curl}\, (\chi \eta ))_i)
    &=
    \epsilon_{ijk} (\delta_{km} -\p_k \eta_m ) \p_j \nabla (\chi \eta_m )
    + 2 \int_0^t\Bigl( \epsilon_{ijk}\p_k v_m \p_j \nabla (\chi
    \eta_m) + \nabla  v (\chi'' \eta + 2  \chi' \nabla \eta )
                \Bigr)   \d s
    \\&\indeq
    + t \chi \nabla \omega_0^i
    +  \nabla \eta (\chi'' \eta + 2\chi' \nabla \eta   ),
  \end{split}
   \llabel{EQ47}
  \end{align}
from where, taking a difference,
  \begin{align}
  \begin{split}
    \nabla ((\curl(\chi E))_i)
    &=
    -
    \epsilon_{ijk} 
    \partial_k E_m \partial_j \nabla (\chi \eta_m)
    +
    \epsilon_{ijk} 
    (\delta_{km} - \partial_k \tilde{\eta}_m) \partial_j \nabla (\chi E_m)
    \\&\indeq
    +
    2 
    \int_0^t \Bigl(     \epsilon_{ijk} \partial_k V_m \partial_j \nabla (\chi \eta_m)
        +    \epsilon_{ijk} \partial_k \tilde{v}_m \partial_j \nabla (\chi E_m)
     \\&\indeq\indeq\indeq\indeq\indeq
	  + \nabla V (\chi'' \eta +2 \chi' \nabla \eta) 
    + \nabla \tilde{v} (\chi'' E + 2\chi' \nabla E)\Bigr)\,\d s 
    \\&\indeq
    + \nabla E (\chi'' \eta + 2 \chi' \nabla \eta) +\nabla \tilde{\eta} (\chi'' E + 2 \chi' \nabla E)
   .
  \end{split}
  \llabel{EQ83}
  \end{align}
Applying the $H^{\delta}$ norm to this identity,  it follows that
\begin{align}
  \begin{split}
  \Vert \nabla \curl(\chi E)\Vert_{{\delta}} 
  &\lec
  \Vert E\Vert_{{1.5+\delta}} 
  + \Vert I- \nabla \tilde{\eta}\Vert_{{1.5+\delta}} \Vert \chi E\Vert_{{2+\delta}}
  + \int_0^t (\Vert V\Vert_{{1.5+\delta}} 
  + \Vert \chi E\Vert_{{2+\delta}} + \Vert E\Vert_{{1+\delta}})\,\d s \\
  &\lec \varepsilon \| \chi E \|_{2+\delta } + \int_0^t Y.
  \end{split}
  \label{EQ84}
\end{align}
The desired estimate \eqref{EQ90} on $\| \chi E\|_{2+\delta }$ now follows directly from applying the div \eqref{EQ82} and curl \eqref{EQ84} estimates, as well as \eqref{EQ85}, in the Brezis-Bourgonion inequality~\eqref{EQ19}.

In order to verify \eqref{EQ96}, we start with the Cauchy invariance
\eqref{EQ42}, which gives 
\begin{align}
  (\curl V)_i = -\epsilon_{ijk}\partial_k E_m \partial_j v_m + \epsilon_{ijk}(\delta_{km}- \partial_k \tilde{\eta}_m)\partial_j V_m
  ,
  \llabel{EQ91}
\end{align}
and the divergence-free condition \eqref{EQ43}  implies that 
\begin{align}
  \div V = -A_{ji} \partial_j v_i + (\delta_{ki} - \tilde{a}_{ji})\partial_j V_i
  .
  \llabel{EQ92}
\end{align}
Consequently,
\begin{align}
  \begin{split}
      &
  \Vert \div V\Vert_{{0.5+\delta}} + \Vert \curl V\Vert_{{0.5+\delta}}
  \lec
  \Vert E\Vert_{{1.5+\delta}} + \Vert A\Vert_{{0.5+\delta}}
  + \Vert V\Vert_{{1.5+\delta}}(\Vert I - \nabla \eta\Vert_{{1.5+\delta}} + \Vert I - \tilde{a}\Vert_{{1.5+\delta}})\\
  &\indeq\lec \varepsilon \| V \|_{1.5+\delta } + \int_0^t \| V \|_{1.5+\delta } \d s,
  \end{split}
  \label{EQ93}
\end{align}
where we used \eqref{EQ85},
$A=\int_{0}^{t}A_t\,\d s$, and
\eqref{EQ87} in the last inequality.
Finally,
\begin{align}
  \begin{split}
  \Vert \nabla_2 V_3\Vert_{H^{\delta}(\Gamma_1)} 
  &\lec \| \Lambda^{0.5} V_3 \|_{H^{0.5+\delta}(\Gamma_1)}
  \lec \| V \|_{L^2} + \varepsilon \Vert V\Vert_{{1.5+\delta}} + \| \Lambda^{0.5} \nabla (\psi V_3) \|_{\delta }
  \\&
  \lec
  \varepsilon \Vert V\Vert_{{1.5+\delta}}
  + \Vert V\Vert_{L^{2}}
  + \Vert \Lambda^{1.5+\delta} (\psi V)\Vert_{L^{2}}
  ,
  \end{split}
  \label{EQ95}
\end{align}
where we used a trace estimate as well as Sobolev interpolation to estimate $\| \Lambda^{0.5} (\psi V_3) \|_{L^2}$ in the first inequality. We also used the divergence-free condition as in \eqref{EQ46} to control $\p_3 V_3$, as well as interpolation again, in the second inequality.

Applying \eqref{EQ93}--\eqref{EQ95} in the Brezis-Bourgonion inequality \eqref{EQ19} gives \eqref{EQ96}, as required.

\subsection{Tangential estimates}\label{sec_uni_tan}

Here we prove~\eqref{EQ118}.
Taking the difference of the Euler equations \eqref{EQ03}
for $v$ and $\tilde v$,
we obtain 
\begin{align}
  \partial_t (\psi V_i) = 
  -A_{ki} \psi \partial_k q
  - \tilde{a}_{ki} \partial_k (\psi Q) + \tilde{a}_{ki} \partial_k \psi Q  
  .
  \llabel{EQ97}
\end{align}
With $S$ defined in \eqref{EQ67}, it follows similarly to \eqref{EQ53} that
\begin{align}
  \begin{split}
    \frac{1}{2}\frac{\d}{\d t}\Vert S(\psi V)\Vert_{L^{2}}^2
    &=
    - \int SA_{ki} \psi \partial_k q S(\psi V_i)\,dx
    - \int \tilde{a}_{ki} \partial_k S(\psi Q) S(\psi V_i)
    \\&\indeq
    + \int S(\tilde{a}_{ki} \partial_k \psi Q) S(\psi V_i)
    - \int (S(A_{ki} \psi \partial_k q) - SA_{ki} \psi \partial_k q) S(\psi V_i)
    \\&\indeq
    - \int (S(\tilde{a}_{ki} \partial_k (\psi Q)) - \tilde{a}_{ki} \partial_k S(\psi Q)) S(\psi V_i)
              \\&
    = I_1 + I_2 + I_3 + I_4 + I_5
    .
  \end{split}
  \label{EQ98}
\end{align}
We start with lower order terms $I_3$--$I_5$. First,
$I_3 \lec \Vert Q\Vert_{{1.5+\delta}} \Vert V\Vert_{{1.5+\delta}}$, while for the commutator terms $I_4$ and $I_5$ we have 
\begin{align}
  \begin{split}
    I_4 &\lec 
    (\Vert \chi^2 A\Vert_{L^{6}}(\Vert \partial_k(\psi q)\Vert_{W^{1.5+\delta,3}} + \Vert \partial_k \psi q\Vert_{W^{1.5+\delta,3}})
             \\&\indeq
    + \Vert \chi^2 A\Vert_{W^{0.5+\delta,3}}(\Vert \partial_k(\psi q)\Vert_{W^{1,6}} + \Vert \partial_k \psi q\Vert_{W^{1,6}}))\Vert V\Vert_{{1.5+\delta}}
    \\&
    \lec
    \Vert \chi^2 A\Vert_{{1+\delta}}\Vert V\Vert_{{1.5+\delta}}
  \end{split}
  \llabel{EQ99}
\end{align}
and 
\begin{align}
  I_5 
  \lec 
  (\Vert \chi^2 \tilde{a}\Vert_{W^{1.5+\delta,3}} \Vert \partial_k(\psi Q)\Vert_{L^{6}} 
  + \Vert \chi^2 \tilde{a}\Vert_{W^{1,6}} \Vert \partial_k(\psi Q)\Vert_{W^{0.5+\delta,3}})\Vert V\Vert_{{1.5+\delta}}
  \lec \Vert \psi Q\Vert_{{2+\delta}}\Vert V\Vert_{{1.5+\delta}}
  .
  \llabel{EQ100}
\end{align}
Thus, $I_3+I_4+I_5$ is bounded by
  \begin{equation}
   C
   (
   \Vert V\Vert_{{1.5+\delta}}
   +
   \Vert \psi Q\Vert_{{2+\delta}}
   +
   \Vert Q\Vert_{{1.5+\delta}}
   +
   \Vert \chi^2 A\Vert_{{1+\delta}}
   +
   \Vert \chi E\Vert_{{2+\delta}}
   )^2
   .
   \label{EQ101}
  \end{equation}
Next, we treat $I_1$ similarly to the fourth term in~\eqref{EQ53}.
We
start by applying $S = \Lambda^{1.5+\delta}$ to the identity
$a\nabla \eta = I$, which leads to
  \begin{align}
  \begin{split}
   Sa \nabla \eta + a S \nabla \eta + r = S I = I
   ,
  \end{split}
   \llabel{EQ102}
  \end{align}
where
$r\coloneqq  S(a \nabla \eta)- (Sa \nabla \eta + a S \nabla \eta)$,
and then
solving for $Sa$ we obtain
  \begin{align}
  \begin{split}
   Sa
    =
    - a S\nabla \eta a
    - r a
    - a
    .
  \end{split}
   \label{EQ103}
  \end{align}
From \eqref{EQ103} we subtract the analogous equation for $\tilde a$,
leading to
  \begin{align}
  \begin{split}
   SA
    =
    - a S\nabla E a
    - A S\nabla \tilde \eta a
    - \tilde a S\nabla \tilde \eta A
    - R a
    - \tilde r A
    - A
    ,
  \end{split}
   \label{EQ104}
  \end{align}
where
$R=r-\tilde r$ and
$\tilde r= S(\tilde a \nabla \tilde\eta)- (S\tilde a \nabla \tilde\eta
+ \tilde a S \nabla \tilde\eta)$.
From \eqref{EQ104}, we then get
\begin{align}
  \begin{split}
   \psi SA
    =
    - a S\nabla(\psi E) a
    + a (S\nabla(\psi E) - \psi S\nabla E  ) a
    - \psi A S\nabla \tilde\eta a
    - \psi \tilde a S\nabla\tilde \eta A
    - \psi R a
    - \psi \tilde r A
    - \psi A
    .
  \end{split}
   \label{EQ105}
  \end{align}
We use \eqref{EQ105} in the expression for $I_1$, which results in
the leading term
  \begin{align}
  \begin{split}
   I_{11}
   &=
    \int
       (a S\nabla (\psi E)a)_{ki}
      \partial_k q S(\psi V_i)
   =
\int
       a_{kl} \partial_{s} S (\psi E_l)a_{si}
      \partial_k q S(\psi V_i)
   \\&
   =
   \int_{\Gamma_1} 
       a_{kl} S (\psi E_l)a_{3i}
      \partial_k q S(\psi V_i)
   -
    \int
       a_{kl} S (\psi E_l)a_{si}
      \partial_k q \partial_{s}S(\psi V_i)
   \\&\indeq
   -
    \int
       a_{kl} S (\psi E_l)a_{si}
      \partial_{sk} q S(\psi V_i)
   -
    \int
       \partial_{s}
       a_{kl} S (\psi E_l)a_{si}
      \partial_k q S(\psi V_i)
   = I_{111} + I_{112} + I_{113} + I_{114}
  \end{split}
   \llabel{EQ106}
  \end{align}
and  seven more terms, whose sum we denote by $\tilde I_{11}$.
It is straight-forward to check that
$\tilde I_{11}$ is bounded by~\eqref{EQ101}.
Similarly, $I_{113}$ and $I_{114}$ are bounded by \eqref{EQ101} and we thus consider them of lower order.
To simplify the presentation, we now introduce the notation indicating equalities and inequalities modulo lower order terms. Thus we write
$I\eqL J$ if $|I-J|$ is bounded from above by~\eqref{EQ101}; 
similarly, $I\leqL J$ if $I-J$ is bounded from above by~\eqref{EQ101}.
With this notation, we may write
$
   I_1
   \eqL
   I_{111} +  I_{112}
$.
For the first term, we have
  \begin{align}
  \begin{split}
   I_{111}
   &=
  \int_{\Gamma_1}
       a_{3m}
       S(\psi E_m)
       a_{3i}
       S(\psi \partial_{t} E_i)
       \partial_{3} q
  \eqL
  \frac12
  \frac{\d}{\d t}
  \int_{\Gamma_1}
       a_{3m}
       S(\psi E_m)
       a_{3i}
       S(\psi E_i)
       \partial_{3} q
   ,
  \end{split}
   \llabel{EQ108}
  \end{align}
while for the second term  $I_{112}$, we use the divergence-free
conditions
$a_{si}\partial_{s} v_i=\tilde a_{si}\partial_{s} \tilde v_i=0$, to
obtain
  \begin{equation}
   a_{si}\partial_{s}V_i
   = - A_{si} \partial_{s}\tilde v_i
  ,
   \llabel{EQ109}
  \end{equation}
whence
\begin{align}
  \begin{split}
   I_{112}
   &\eqL
   -  \int a_{km}
       S(\psi E_m)
       \partial_{k} q
      S(a_{si}\psi \partial_{s}V_i)
   \eqL
     \int a_{km}
       S(\psi E_m)
       \partial_{k} q
      S(\chi^2 A_{si}\psi \partial_{s}\tilde v_i)
   \\&
   =
     \int
      \Lambda^{1/2}
      (
       a_{km}
       S(\psi E_m)
       \partial_{k} q
       )
      \Lambda^{-1/2}
      S(\chi^2 A_{si}\psi \partial_{s}\tilde v_i)
   \eqL 0
   .
  \end{split}
   \llabel{EQ110}
  \end{align}
Using the above estimates, we obtain
  \begin{align}
  \begin{split}
   I_1
    &\eqL
      \frac12
  \frac{\d}{\d t}
  \int_{\Gamma_1}
       a_{3m}
       S(\psi E_m)
       a_{3i}
       S(\psi E_i)
       \partial_{k} q
    .
  \end{split}
   \llabel{EQ111}
  \end{align}

In order to bound $I_2 = - \int \tilde{a}_{ki} \partial_k S(\psi Q) S(\psi V_i)$,
we integrate by parts in $x_k$ to obtain 
\begin{align}
  \begin{split}
  I_2 &= \int S(\psi Q) \tilde{a}_{ki} \partial_k S(\psi V_i) 
  \eqL
      \int S(\psi Q) S(\tilde{a}_{ki}\psi \partial_{k}V_i) 
        \\&=\int S(\psi Q) S(\chi^2{A}_{ki}\psi \partial_{k}v_i) 
        =  \int \Lambda^{1/2}S(\psi Q)  \Lambda^{-1/2}S(\chi^2{A}_{ki}\psi \partial_{k}v_i) 
  \eqL0
  .
  \end{split}
  \llabel{EQ112}
\end{align}
Now, we integrate \eqref{EQ98} in time and collect all the tangential estimates
and bound
\begin{align}
  \Vert S(\psi V)\Vert_{L^{2}}^2+\Vert a_{3i}SE_i\Vert_{L^{2}(\Gamma_1)}^2
  \lec
  \int_0^t
  \Bigl(Y^2
  +\Vert \psi Q\Vert_{{2+\delta}}^2 +\Vert Q\Vert_{{1.5+\delta}}^2
  \Bigr)\,\d s
  ,
  \label{EQ113}
\end{align}
where we used the Rayleigh-Taylor condition \eqref{EQ15_RT} to obtain the second term on the left-hand side. 

In order to estimate the pressure terms in \eqref{EQ113}, we first note that $Q$ satisfies the Poisson equation 
\begin{align}
  \begin{split}
    \Delta  Q 
    &=
    -\partial_j((A_{ji}a_{ki} + \tilde{a}_{ji}A_{ki})\partial_k q)
    + \partial_j((\delta_{jk}-\tilde{a}_{ji}\tilde{a}_{jk})\partial_k Q)
    \\&\indeq
    +\p_j ( \partial_t A_{ji}  v_i )
    + \p_j (\partial_t \tilde{a}_{ji} V_i )
    \llabel{EQ114a}   
  \end{split}  
\end{align}
(see \eqref{EQ27}),
where we also used the Piola identity $\p_j A_{ji}=0$. The elliptic regularity in $H^{1.5+\delta }$ \cite[Chapter~2]{LM} then gives 
\eqnb\llabel{EQ114b}
\begin{split}
\| Q \|_{1.5+\delta } &\lec \| (Aa + \tilde a A) \na q\|_{0.5+\delta } + \| (I-\tilde a \tilde a^T)\na Q \|_{0.5+\delta } + \| A_t v \|_{0.5+\delta } + \| \tilde a_t V \|_{0.5+\delta }\\
&\lec \| A \|_{0.5+\delta } + \varepsilon \| Q \|_{1.5+\delta } + \| V \|_{1.5+\delta } + \| E \|_{1.5+\delta } 
,
\end{split}
\eqne
where we used \eqref{EQ87} in the second inequality. 
Absorbing $\varepsilon \| Q \|_{1.5+\delta }$ by the left-hand side and using  \eqref{EQ85}--\eqref{EQ86}  we thus obtain 
\eqnb\label{EQ117a}
\| Q \|_{1.5+\delta } \lec Y+\int_0^t Y.
\eqne

In order to estimate $\| \psi Q \|_{2+\delta}$ we note that $\psi Q$ solves the Poisson equation 
\begin{align}
  \begin{split}
    \Delta (\psi Q) 
    &=
    -\partial_j((A_{ji}a_{ki}+\tilde{a}_{ji}A_{ki})\partial_k(\psi q))
    + \partial_j((\delta_{jk}-\tilde{a}_{ji}\tilde{a}_{jk})\partial_k(\psi Q))
    \\&\indeq
    +\partial_j(A_{ji}a_{ki}\partial_k \psi q)
    +\partial_j(\tilde{a}_{ji}A_{ki}\partial_k\psi q)
    + \partial_j(\tilde{a}_{ji}\tilde{a}_{ki}\partial_k\psi Q)
    \\&\indeq
    +\psi \partial_t A_{ji} \partial_j v_i
    +\psi \partial_t \tilde{a}_{ji}\partial_j V_i
    +A_{ji}\partial_j \psi a_{ki}\partial_k q 
    +\tilde{a}_{ji}\partial_j \psi A_{ki}\partial_k q 
    +\tilde{a}_{ji}\partial_j \psi \tilde{a}_{ki} \partial_k Q  
    ,
    \label{EQ114}   
  \end{split}  
\end{align}
with homogeneous Dirichlet boundary conditions,  see \cite[(3.18)]{KO}
for details.
Therefore, by the elliptic regularity, we obtain
\begin{align}
  \begin{split}
    \Vert \psi Q\Vert_{{2+\delta}}
    &\lec
    \Vert Aa\nabla(\psi q)\Vert_{{1+\delta}}
    +\Vert \tilde a A\nabla(\psi q)\Vert_{{1+\delta}}
    +\Vert (I-\tilde{a}\tilde{a}^T)\nabla (\psi Q)\Vert_{{1+\delta}}
    \\&\indeq
    +\Vert Aa\nabla\psi q\Vert_{{1+\delta}}
    +\Vert \tilde{a}A\nabla\psi q\Vert_{{1+\delta}}
    +\Vert \tilde{a}^2\nabla \psi Q\Vert_{{1+\delta}}
    \\&\indeq
    +\Vert \psi A_t \nabla v\Vert_{{\delta}}
    +\Vert \psi \tilde{a}_t\nabla V\Vert_{{\delta}}
    +\Vert A\nabla \psi a \nabla q\Vert_{{\delta}}
    +\Vert \tilde{a}\nabla \psi A\nabla q\Vert_{{\delta}}
    +\Vert \tilde a^2\nabla \psi \nabla Q\Vert_{{\delta}}
    ,  
    \llabel{EQ115}
  \end{split}
\end{align}
whence
\begin{align}
  \begin{split}
    &
    \Vert \psi Q\Vert_{{2+\delta}}\lec
    \Vert \chi^2 A\Vert_{{1+\delta}}
    +\Vert I-\tilde{a}\tilde{a}^T\Vert_{{1.5+\delta}}\Vert \psi Q\Vert_{{2+\delta}}
    +\Vert Q\Vert_{{1+\delta}}
    +\Vert A_t\Vert_{{\delta}}+\Vert V\Vert_{{1+\delta}} 
    \\&\indeq
    \lec
    \Vert \chi E\Vert_{{2+\delta}}+\Vert Q\Vert_{{1+\delta}}+\Vert V\Vert_{{1+\delta}}+ \varepsilon \| \psi Q \|_{2+\delta}
    ,
    \llabel{EQ115}
  \end{split}
\end{align}
where we used \eqref{EQ86} to bound $\| \chi^2 A \|_{1+\delta }$ and
\eqref{EQ87} to estimate $\| \p_t A \|_{\delta }$.
Absorbing the last term by the left-hand side and applying \eqref{EQ117a} we obtain
\[
\| \psi Q \|_{2+\delta } \lec Y+ \int_0^t Y
.
\]
Applying this and \eqref{EQ117a} in \eqref{EQ113} gives 
\begin{align}
  \Vert S(\psi V)\Vert_{L^{2}}^2 + \Vert a_{3i}SE_i\Vert^2_{L^{2}(\Gamma_1)}
  \lec 
  \int_0^t Y^2
  \label{EQ118}
  .
\end{align}
Finally, writing
\begin{align}
  \Vert  SE_3\Vert_{L^{2}(\Gamma_1)}
  \lec 
  \Vert a_{3i}SE_3\Vert_{L^{2}(\Gamma_1)} + \Vert I-a\Vert_{{1.5+\delta}}\Vert \chi E\Vert_{{2+\delta}}   \lec \Vert a_{3i}SE_3\Vert_{L^{2}(\Gamma_1)} + \varepsilon \Vert \chi E\Vert_{{2+\delta}}  
  ,
  \llabel{EQ120}
\end{align}
we obtain \eqref{EQ118}, as required.

\subsection{$L^2$ control of $V$}\label{sec_uni_l2}
Here we show~\eqref{EQ121}. 

To this end, we note that the difference of the Euler equations \eqref{EQ03} gives $V_t = -A \nabla q - \tilde{a} \nabla Q$, and taking the $L^2$ norm gives 
\[ \Vert V_t \Vert_{L^{2}} \lec \Vert A\Vert_{L^{2}} + \Vert Q\Vert_{{1}} \lec \| E \|_{1+\delta } + \int_0^t Y \lec \int_0^t Y , \]
as required, where we used \eqref{EQ86} and~\eqref{EQ117a}.

\section*{Acknowledgments} 
MSA and IK were supported in part by the NSF grant DMS-2205493, while
WSO was supported in part by the Simons Foundation.


\begin{thebibliography}{[WWZ]}
\small
\bibitem[ABZ1]{ABZ1} 
T.~Alazard, N.~Burq, and C.~Zuily, \emph{On the water-wave equations
with surface tension}, Duke Math.~J.~\textbf{158} (2011), no.~3,
413--499.

\bibitem[ABZ2]{ABZ2} 
T.~Alazard, N.~Burq, and C.~Zuily, \emph{Low regularity {C}auchy theory for the
  water-waves problem: canals and wave pools}, Lectures on the analysis of
  nonlinear partial differential equations. {P}art 3, Morningside Lect. Math.,
  vol.~3, Int. Press, Somerville, MA, 2013, pp.~1--42. 
 
 
\bibitem[AD]{AD} 
T.~Alazard and J.-M.~Delort, \emph{Global solutions and asymptotic
  behavior for two dimensional gravity water waves}, Ann. Sci. \'Ec. Norm.
  Sup\'er. (4) \textbf{48} (2015), no.~5, 1149--1238. 

\bibitem[AM1]{AM1} 
D.M.~Ambrose and N.~Masmoudi, \emph{The zero surface tension limit of
  two-dimensional water waves}, Comm. Pure Appl. Math. \textbf{58} (2005),
  no.~10, 1287--1315. 

\bibitem[AM2]{AM2} 
D.M.~Ambrose and N.~Masmoudi, \emph{The zero surface tension limit of
  three-dimensional water waves}, Indiana Univ. Math. J. \textbf{58} (2009),
  no.~2, 479--521. 

\bibitem[B]{B} 
J.T.~Beale, \emph{The initial value problem for the {N}avier-{S}tokes
  equations with a free surface}, Comm. Pure Appl. Math. \textbf{34} (1981),
  no.~3, 359--392. 
  


\bibitem[BL1]{BL1} 
J.~Bourgain and D.~Li,
  \emph{Strong ill-posedness of the incompressible Euler equation in borderline Sobolev spaces}, Invent. math. \textbf{201} (2015), 97--157.

\bibitem[BL2]{BL2} 
J.~Bourgain and D.~Li,
  \emph{Strong ill-posedness of the 3D incompressible Euler equation in borderline spaces}, Int. Math. Res. Not. \textbf{16} (2021), 12155--12264.


\bibitem[BB]{EQ19} 
J.P.~Bourguignon and H.~Brezis,
  \emph{Remarks on the {E}uler equation}, J.~Functional Analysis \textbf{15} (1974), 341--363.
  

\bibitem[CL]{CL} 
A.~Castro and D.~Lannes, \emph{Well-posedness and shallow-water stability
  for a new {H}amiltonian formulation of the water waves equations with
  vorticity}, Indiana Univ. Math. J. \textbf{64} (2015), no.~4, 1169--1270.

\bibitem[ChL]{ChL} 
D.~Christodoulou and H.~Lindblad, \emph{On the motion of the free
  surface of a liquid}, Comm. Pure Appl. Math. \textbf{53} (2000), no.~12,
  1536--1602. 

\bibitem[CS1]{CS1} 
D.~Coutand and S.~Shkoller, \emph{Well-posedness of the free-surface
  incompressible {E}uler equations with or without surface tension}, J. Amer.
  Math. Soc. \textbf{20} (2007), no.~3, 829--930.

\bibitem[CS2]{CS2} 
D.~Coutand and S.~Shkoller, \emph{A simple proof of well-posedness for
  the free-surface incompressible {E}uler equations}, Discrete Contin. Dyn.
  Syst. Ser. S \textbf{3} (2010), no.~3, 429--449. 
  
\bibitem[DE1]{DE1} 
M.M.~Disconzi and D.G.~Ebin, 
  \emph{The free boundary {E}uler equations with large surface tension}, 
  J. Differential Equations \textbf{261} (2016), no.~2, 821--889. 

\bibitem[DE2]{DE2} 
M.M.~Disconzi and D.G.~Ebin, 
  \emph{On the limit of large surface tension for a fluid motion with free boundary}, 
  Comm. Partial Differential Equations \textbf{39} (2014), no.~4, 740--779. 

\bibitem[DK]{DK}
  M.M.~Disconzi and I.~Kukavica,
  \emph{A priori estimates for the free-boundary Euler equations with surface tension in three dimensions}
  (submitted).

\bibitem[E]{E} 
D.G.~Ebin, \emph{The equations of motion of a perfect fluid with free
  boundary are not well posed}, Comm. Partial Differential Equations
  \textbf{12} (1987), no.~10, 1175--1201. 
  
\bibitem[EL]{EL}
T.~Elgindi and D.~Lee, \emph{Uniform regularity for free-boundary Navier-Stokes equations with surface tension},
arXiv:1403.0980 (2014).

  
\bibitem[GMS]{GMS} 
P.~Germain, N.~Masmoudi, and J.~Shatah, \emph{Global solutions for the gravity
  water waves equation in dimension 3}, Ann. of Math. (2) \textbf{175} (2012),
  no.~2, 691--754. 

\bibitem[HIT]{HIT} 
J.K.~Hunter, M.~Ifrim, and D.~Tataru, \emph{Two dimensional water
  waves in holomorphic coordinates}, Comm. Math. Phys. \textbf{346} (2016),
  no.~2, 483--552. 

\bibitem[IT]{IT} 
M.~Ifrim and D.~Tataru, \emph{Two dimensional water waves in
  holomorphic coordinates {II}: {G}lobal solutions}, Bull. Soc. Math. France
  \textbf{144} (2016), no.~2, 369--394. 

\bibitem[I]{I} 
T.~Iguchi, \emph{Well-posedness of the initial value problem for
  capillary-gravity waves}, Funkcial. Ekvac. \textbf{44} (2001), no.~2,
  219--241. 

\bibitem[IK]{IK} 
M. Ignatova and I. Kukavica, \emph{On the local existence of the
  free-surface {E}uler equation with surface tension}, Asymptot. Anal.
  \textbf{100} (2016), no.~1-2, 63--86. 
  
\bibitem[IP]{IP} 
A.D.~Ionescu and F.~Pusateri, \emph{Global solutions for the gravity
  water waves system in 2d}, Invent. Math. \textbf{199} (2015), no.~3,
  653--804. 

\bibitem[KO]{KO}
I.~Kukavica and W.S.~O\.za\'nski,
  \emph{Local-in-time existence of free-surface 3D Euler flow with $H^{2+\delta}$ initial vorticity in a neighborhood of the free boundary}, Nonlinearity \textbf{36} (2023), no.~1, 636--652.


\bibitem[KP]{KP} 
T.~Kato and G.~Ponce, \emph{Commutator estimates and the {E}uler and
  {N}avier-{S}tokes equations}, Comm. Pure Appl. Math. \textbf{41} (1988),
  no.~7, 891--907. 

\bibitem[KT1]{KT1} 
I.~Kukavica and A.~Tuffaha, \emph{On the 2{D} free boundary {E}uler
  equation}, Evol. Equ. Control Theory \textbf{1} (2012), no.~2, 297--314.

\bibitem[KT2]{KT2} 
I.~Kukavica and A.~Tuffaha, \emph{A regularity result for the
  incompressible {E}uler equation with a free interface}, Appl. Math. Optim.
  \textbf{69} (2014), no.~3, 337--358. 

\bibitem[KT3]{KT3} 
I.~Kukavica and A.~Tuffaha, \emph{Well-posedness for the compressible
  {N}avier-{S}tokes-{L}am\'e system with a free interface}, Nonlinearity
  \textbf{25} (2012), no.~11, 3111--3137. 

\bibitem[KT4]{KT4} 
I.~Kukavica and A.~Tuffaha, 
  \emph{A sharp regularity result for the {E}uler equation with a free interface}, 
  Asymptot. Anal. \textbf{106} (2018), no.~2, 121--145. 

\bibitem[KTV]{KTV}
I.~Kukavica, A.~Tuffaha, and V.~Vicol,
{\em On the local existence and uniqueness for the 3D Euler equation
with a free interface},
 Appl.\ Math.\ Optim. (2016), doi:10.1007/s00245-016-9360-6.

\bibitem[KTVW]{KTVW} 
I.~Kukavica, A.~Tuffaha, V.~Vicol, and F.~Wang, \emph{On the existence
  for the free interface 2{D} {E}uler equation with a localized vorticity
  condition}, Appl. Math. Optim. \textbf{73} (2016), no.~3, 523--544.

\bibitem[L]{L} 
D.~Lannes, \emph{Well-posedness of the water-waves equations}, J. Amer.
  Math. Soc. \textbf{18} (2005), no.~3, 605--654 (electronic). 

\bibitem[Li]{Li} 
D.~Li, 
  \emph{On {K}ato-{P}once and fractional {L}eibniz}, 
  Rev.\ Mat.\ Iberoam.~\textbf{35} (2019), no.~1, 23--100. 
  
\bibitem[Lin1]{Lin1} 
H.~Lindblad, \emph{Well-posedness for the linearized motion of an
  incompressible liquid with free surface boundary}, Comm. Pure Appl. Math.
  \textbf{56} (2003), no.~2, 153--197.

\bibitem[Lin2]{Lin2} 
H.~Lindblad, \emph{Well-posedness for the motion of an incompressible liquid
  with free surface boundary}, Ann. of Math. (2) \textbf{162} (2005), no.~1,
  109--194. 

\bibitem[LM]{LM} 
  J.-L.~Lions and E.~Magenes, 
  \emph{Nonhomogeneous boundary value problems and applications: vol. 1}, 
   Springer-Verlag Berlin Heidelberg, 1972.

\bibitem[MC]{MC} 
B.~Muha and S.~\v Cani\'c, \emph{Fluid-structure interaction between
  an incompressible, viscous 3{D} fluid and an elastic shell with nonlinear
  {K}oiter membrane energy}, Interfaces Free Bound. \textbf{17} (2015), no.~4,
  465--495. 

  
\bibitem[MR]{MR} 
N.~Masmoudi and F.~Rousset, \emph{Uniform regularity and vanishing
  viscosity limit for the free surface {N}avier-{S}tokes equations}, Arch.
  Ration. Mech. Anal. \textbf{223} (2017), no.~1, 301--417. 


\bibitem[N]{N} 
V.I.~Nalimov, \emph{The Cauchy-Poisson problem}, Dinamika Splo\v sn. Sredy
  (1974), no.~Vyp. 18 Dinamika Zidkost. so Svobod. Granicami, 104--210, 254.

\bibitem[OT]{OT} 
M.~Ogawa and A.~Tani, \emph{Free boundary problem for an incompressible
  ideal fluid with surface tension}, Math. Models Methods Appl. Sci.
  \textbf{12} (2002), no.~12, 1725--1740. 

\bibitem[P]{P} 
F.~Pusateri, \emph{On the limit as the surface tension and density ratio
  tend to zero for the two-phase {E}uler equations}, J. Hyperbolic Differ. Equ.
  \textbf{8} (2011), no.~2, 347--373. 

\bibitem[S]{S} 
B.~Schweizer, \emph{On the three-dimensional {E}uler equations with a free
  boundary subject to surface tension}, Ann. Inst. H. Poincar\'e Anal. Non
  Lin\'eaire \textbf{22} (2005), no.~6, 753--781.

\bibitem[Sh]{Sh} 
M.~Shinbrot, \emph{The initial value problem for surface waves under
  gravity. {I}. {T}he simplest case}, Indiana Univ. Math. J. \textbf{25}
  (1976), no.~3, 281--300. 

\bibitem[Shn]{Shn} 
A.I.~Shnirelman, \emph{The geometry of the group of diffeomorphisms and the
  dynamics of an ideal incompressible fluid}, Mat. Sb. (N.S.) \textbf{128(170)}
  (1985), no.~1, 82--109, 144. 

\bibitem[SZ1]{SZ1} 
J.~Shatah and C.~Zeng, \emph{Geometry and a priori estimates for free
  boundary problems of the {E}uler equation}, Comm. Pure Appl. Math.
  \textbf{61} (2008), no.~5, 698--744. 

\bibitem[SZ2]{SZ2} 
J.~Shatah and C.~Zeng, \emph{Local well-posedness for fluid interface
  problems}, Arch. Ration. Mech. Anal. \textbf{199} (2011), no.~2, 653--705.

\bibitem[T]{T} 
A.~Tani, \emph{Small-time existence for the three-dimensional
  Navier-Stokes equations for an incompressible fluid with a free surface},
  Arch. Rational Mech. Anal. \textbf{133} (1996), no.~4, 299--331. 


\bibitem[WZZZ]{WZZZ}
C.~Wang, Z.~Zhang, W.~Zhao, and Y.~Zheng, \emph{Local well-posedness and break-down criterion of the incompressible Euler equations with free boundary},   arXiv:1507.02478, 2015.

\bibitem[W1]{W1} 
S.~Wu, \emph{Well-posedness in Sobolev spaces of the full water wave
  problem in {$2$}-{D}}, Invent. Math. \textbf{130} (1997), no.~1, 39--72.

\bibitem[W2]{W2} 
S.~Wu, \emph{Well-posedness in Sobolev spaces of the full water wave
  problem in 3-{D}}, J. Amer. Math. Soc. \textbf{12} (1999), no.~2, 445--495.
  
\bibitem[W3]{W3} 
S.~Wu, \emph{Global wellposedness of the 3-{D} full water wave problem},
  Invent. Math. \textbf{184} (2011), no.~1, 125--220. 


\bibitem[Y1]{Y1} 
H.~Yosihara, \emph{Gravity waves on the free surface of an incompressible
  perfect fluid of finite depth}, Publ. Res. Inst. Math. Sci. \textbf{18}
  (1982), no.~1, 49--96. 

\bibitem[Y2]{Y2} 
H.~Yosihara, \emph{Capillary-gravity waves for an incompressible ideal
  fluid}, J. Math. Kyoto Univ. \textbf{23} (1983), no.~4, 649--694. 
  

\bibitem[ZZ]{ZZ} 
P.~Zhang and Z.~Zhang, \emph{On the free boundary problem of
  three-dimensional incompressible {E}uler equations}, Comm. Pure Appl. Math.
  \textbf{61} (2008), no.~7, 877--940. 


\end{thebibliography}
\end{document}